\documentclass[12pt, a4paper, twoside]{article}
\usepackage{amsthm, amsfonts, amsmath, amscd, amssymb, epsfig}
\usepackage[left=2.5cm, right=2.5cm]{geometry}
\usepackage[UKenglish]{babel}
\usepackage[utf8]{inputenc}
\usepackage{hyperref}
\usepackage{german}
\usepackage{epigraph} 
\numberwithin{equation}{section}
\theoremstyle{plain}
\newtheorem{theorem}{Theorem}[section]

\newtheorem{lemma}[theorem]{Lemma}
\newtheorem{lemma-definition}[theorem]{Lemma-Definition}

\theoremstyle{definition}
\newtheorem{definition}[theorem]{Definition}

\newtheorem{remark}[theorem]{Remark}

\newtheorem{problem}[theorem]{Problem}

\newcommand{\Addresses}{{
  \bigskip
\footnotesize
\textsc{Shah Faisal\\ Department of Mathematics
\\
Uppsala University\\ Sweden}\par\nopagebreak
\textit{E-mail address:} \texttt{shahmath19@gmail.com,  shah.faisal@math.uu.se}}}
\graphicspath{{Images/}}

\begin{document}

\title{Stabilized symplectic embeddings of higher-dimensional ellipsoids}
\author{Shah Faisal}
\date{}
\maketitle
\begin{abstract}
\noindent We provide a lower bound for the embedding capacity of higher-dimensional symplectic ellipsoids, formulated in terms of the Lagrangian capacity of ellipsoids. Our approach relies on examining the Borman--Sheridan class of a Weinstein neighborhood of a suitable monotone Lagrangian torus, using Tonkonog's string topology-based computation of the gravitational descendants of the torus. \end{abstract}	
\parskip=4pt
\section{Introduction and statements of results}
A symplectic embedding problem asks whether one symplectic manifold can be symplectically embedded into another. Beyond their intrinsic interest, symplectic embeddings can be used to define new global invariants for symplectic manifolds \cite{MR3777016, MR2369441}. Due to the co-existence of rigidity and flexibility in symplectic embeddings, studying symplectic embedding problems helps to illustrate the boundary between rigidity and flexibility in symplectic geometry.  Symplectic embedding techniques have also found intriguing applications to classical dynamical systems, such as the restricted three-body problem \cite{Frauenfelder_2018}; see also {\cite[Section 4.5]{MR3777016}}.

In this article, we consider the problem of embedding one standard ellipsoid into another. The standard $2n$-dimensional ellipsoid is defined by
\[E^{2n}(x_1,x_2,\dots,x_n):=\bigg\{(z_1,\dots, z_{n})\in \mathbb{C}^{n}: \sum_{i=1}^{n}\frac{  \pi|z_i|^2}{x_i}\leq 1\bigg \},\] 
where $0<x_1\leq x_2\leq x_3\leq \dots\leq x_n\leq\infty$. In particular,
\[\bar{B}^{2n}(r):=E^{2n}(r,r,\dots,r)\] 
is the standard closed ball of capacity $r>0$ and radius $\sqrt{r/\pi}$. Moreover, $S^{2n-1}(r):=\partial \bar{B}^{2n}(r)$ is the sphere of radius $\sqrt{r/\pi}$ centered at the origin.  In particular, $S^{1}(r)$ is the circle centered at the origin that bounds a symplectic area equal to $r$. We denote by $B^{2n}(r)$ the open ball of capacity $r$.

 Ellipsoids inherit the standard symplectic structure from $(\mathbb{C}^{n}, \omega_{\mathrm{std}}:=\frac{i}{2}\sum_{j=1}^ndz_j\wedge d\bar{z}_j)$. This way, they form an essential class of symplectic manifolds with boundary. 
 
 The symplectic embedding problem for the standard ellipsoids is stated as follows.
\begin{problem}\label{symellipproblem}
Given $0<x_1\leq x_2\leq \dots\leq x_n\leq\infty$, determine all $0<b_1\leq b_2\leq \dots\leq b_n\leq\infty$ such that $E^{2n}(x_1,x_2,\dots,x_n)$ embeds into  $E^{2n}(b_1,b_2,\dots,b_n)$ symplectically.  
\end{problem}
Problem \ref{symellipproblem} has been completely solved in dimension four \cite{McDuff-Hoferconjecture, McDuff-Schlenk,Faisal:2024aa}. However, little is known in higher dimensions; we recommend \cite{Buse-Hindellipsoidsindimensiongreaterthanfour, Cristofaro-Hind-Siegel-stabilizedembedding, MR3777016,  Schlenkellipsoids} for some progress and open problems.

In this paper, we consider the following two cases of  Problem \ref{symellipproblem}.

\textbf{Case 1:} Let $n\geq 2$ be a positive integer. In this case, we are interested in determining the embedding capacity defined by
\begin{equation}\label{function-embed0}
\mathcal{EB} _{n}(x_2,x_3,\dots,x_{n}):=\inf\big\{r>0: E^{2n}(1, x_2, \dots, x_{n})\xrightarrow[]{s} B^{2n}(r)\big\}
\end{equation}
for any $1\leq x_2\leq x_3\leq \dots \leq x_n< \infty$, 
where ``$\xrightarrow[]{s}$'' represents a symplectic embedding with the standard symplectic structures on the domain and the target. 

\textbf{Case 2:}   Let $n\geq 3$ and $k\geq 2$ be positive integers such that $n-k\geq 1$. In this case, the problem asks to determine the embedding capacity defined by
\begin{equation}\label{function-embed}
\mathcal{EC} _{n-k}(x_2,x_3,\dots,x_{n}):=\inf\big\{r>0: E^{2n}(1, x_2, \dots, x_{n})\xrightarrow[]{s} B^{2k}(r)\times \mathbb{C}^{n-k}\big\}
\end{equation}
for any $1\leq x_2\leq x_3\leq \dots \leq x_n\leq \infty$, 
where ``$\xrightarrow[]{s}$'' represents a symplectic embedding with the standard symplectic structures on both the domain and target.  We mention that  function (\ref{function-embed}), for $k=2$, has already appeared in the work of Hind {\cite[Section 1.2]{MR3394319}}.

The (restricted) stabilized embedding problem for four-dimensional ellipsoids  is about determining (\ref{function-embed}) in the special situation: \[k=2 \text{ and } x_3=\dots=x_n=\infty,\]

 i.e., the function $\mathcal{EC} _{n-2}(x_2,\infty,\dots,\infty):[1,\infty]\to \mathbb{R}_{\geq 0}$ defined by
\begin{equation}\label{emb-fuction00}
\mathcal{EC} _{n-2}(x_2,\infty,\dots,\infty):=\inf\big\{r>0: E^4(1,x)\times \mathbb{C}^{n-2}\xrightarrow[]{s} B^{4}(r)\times \mathbb{C}^{n-2} \big\}
\end{equation}
for any  integer $n\geq 3$.

We briefly mention some notable progress on the above two cases that we are aware of.
\begin{itemize}
\item McDuff--Schlenk \cite{McDuff-Schlenk}: the function $\mathcal{EB} _{2}(\cdot)$ is understood completely. On the interval $[0,\tau^4]$, where $\tau=\frac{1+\sqrt{5}}{2}$, the function $\mathcal{EB} _{2}(x_2)$ is given by the Fibonacci stairs. For any $x_2\in [\tau^4, 7]$, we have $3\mathcal{EB} _{2}(x_2)=(x_2+1)$. For all $x_2\geq 8\frac{1}{36}$, we have $\mathcal{EB} _{2}(x_2)=\sqrt{x_2}$. On the interval $[7,8\frac{1}{36}]$, the function $\mathcal{EB} _{2}(x_2)$ takes a complicated form; see {\cite[Theorem 1.1.2]{McDuff-Schlenk}}. 
 \item Buse--Hind \cite{Buse-Hindellipsoidsindimensiongreaterthanfour}: for six-dimensional ellipsoids $E^{6}(1, x_2, x_3)$ satisfying  $x_2^2+x_3^2\geq 1.41\times 10^{101}$, the volume constraint is optimal: $\mathcal{EB} _{3}(x_2, x_3)=\sqrt[3]{x_2x_3}$. Also for small values of $x_2$ and $x_3$, the values $\mathcal{EB} _{3}(x_2, x_3)$ are known; see {\cite[Lemma 2.7---2.11, Figure 1]{Buse-Hindellipsoidsindimensiongreaterthanfour}}.

\item Cristofaro--Hind--McDuff \cite{MR3789827}: on the interval $[0,\tau^4]$, where $\tau=\frac{1+\sqrt{5}}{2}$, the function $\mathcal{EC} _{n-2}(x_2,\infty,\dots,\infty)$ is given by the Fibonacci stairs for any $n\geq 3$.
\item Hind \cite{MR3394319}: For $x_2>\tau^4$ and any $n\geq 3$ it holds that
\begin{equation}\label{hind}
	\mathcal{EC} _{n-2}(x_2,\infty,\dots,\infty)\leq \frac{3x_2}{x_2+1}.
\end{equation} 
\item  Hind--Kerman \cite{Hind-Kerman-ellipembedding,MR3867635}: for $x_2=3l-1$, where $l$ is an odd indexed Fibonacci number, we have 
\[\mathcal{EC} _{n-2}(x_2,\infty,\dots,\infty)=\frac{3x_2}{x_2+1}\]
 for any $n\geq 3$.
\item McDuff \cite{MR3782228}: for $x_2=3l-1$, where $l$ is any positive integer, we have
 \[\mathcal{EC} _{n-2}(x_2,\infty,\dots,\infty)=\frac{3x_2}{x_2+1}\]
 for any $n\geq 3$.
 \item McDuff--Siegel \cite{McDuff:2024aa}:  for any integer $n\geq 3$ and all $x_2\geq\tau^4$, we have 
\[ \mathcal{EC} _{n-2}(x_2,\infty,\dots,\infty)\geq \frac{3x_2}{x_2+1}.\] 
\end{itemize}

From the above works, we conclude that the function $\mathcal{EC} _{n-2}(x_2,\infty,\dots,\infty)$ is now fully understood for any $n\geq 3$.  On the interval $[0,\tau^4]$, where $\tau=\frac{1+\sqrt{5}}{2}$, the function $\mathcal{EC} _{n-2}(x_2,\infty,\dots,\infty)$ is given by the Fibonacci stairs for any $n\geq 3$. For all $n\geq 3$ and all $x_2\geq\tau^4$ we have
  \begin{equation}\label{SE4}
 \mathcal{EC} _{n-2}(x_2,\infty,\dots,\infty)= \frac{3x_2}{x_2+1}.
 \end{equation}
 \subsection{Main results}
We now state our main theorems. 
\begin{theorem}\label{maintheorem2}
For any positive integers $n$ and $k$ such that $n-k\geq 1$ and $k>2$, there exists $0<c_k<\infty$ such that
\begin{equation}\label{SEh01}
\mathcal{EC} _{n-k}(x_2,x_3,\dots,x_{n})\geq (k+1)\bigg(1+\frac{1}{x_2}+\cdots+\frac{1}{x_n}\bigg)^{-1}.
\end{equation}
for all  $c_k\leq x_2\leq x_3\leq \dots\leq x_n\leq \infty$. Moreover, we have $\lim_{k\to \infty}c_k=\infty$.
\end{theorem}
Inspired by Gromov's pioneering idea in the proof of his non-squeezing theorem, one can seek sharp numerical obstructions to symplectic embeddings by constructing pseudo-holomorphic curves with carefully controlled symplectic area. In case of stabilized embeddings of ellipsoids, this approach leads to the problem of producing pseudo-holomorphic curves of arbitrarily large degree, subject to specific geometric and asymptotic constraints, in complex projective spaces; see \cite{Hind-Kerman-ellipembedding, McDuff:2024aa, Mikhalkin:2023aa, MR4466005, McDuff:2023aa} and the references therein. In particular, a proof of (\ref{SE4}) reduces to establishing the existence of certain curves of arbitrarily large degree carrying a special type of singularities   {\cite{MR3782228,  McDuff:2024ab, MR4466005, McDuff:2024aa}}.

Constructing pseudo-holomorphic curves with a prescribed symplectic area is often highly challenging because of the area constraint. To avoid this difficulty, we adopt a different point of view: rather than searching for curves of a fixed area, we observe that the very existence of a symplectic embedding typically forces the existence of certain rigid pseudo-holomorphic curves that can be counted (cf. Theorem \ref{keytheorem1} and Theorem \ref{importantcount}). These curve counts can therefore be regarded as obstructions to the existence of particular symplectic embeddings: an embedding that surpasses the expected obstruction would necessarily imply an impossible count of curves--either too many or too few. In this sense, the obstruction arises from curve counts rather than symplectic areas. This is the guiding idea that we follow in the proof of Theorem \ref{maintheorem2}.
\begin{remark}
We note that Theorem \ref{maintheorem2} does not cover the case $k=2$, but the statement holds in this case for $x_3=\cdots=x_n=\infty$ and $n\geq 3$.  In fact, (\ref{SE4}) implies that (\ref{SEh01}) becomes an equality for $k=2$, $x_3=\cdots=x_n=\infty$ and $x_2\geq\tau^4$:
\[\mathcal{EC} _{n-2}(x_2,\infty,\dots,\infty)=\frac{3x_2}{x_2+1}=3\bigg(1+\frac{1}{x_2}\bigg)^{-1}.\]
\end{remark}
\begin{remark}
By {\cite[Theorem 4.37]{Pereira:2025aa}}, the Lagrangian capacity of an ellipsoid is given by
\begin{equation}\label{Lagcapellip}
\operatorname{C}_{\mathrm{Lag}}(E^{2n}(1,x_2,x_3,\dots,x_{n}),\omega_{\mathrm{std}})=\bigg(1+\frac{1}{x_2}+\cdots+\frac{1}{x_n}\bigg)^{-1}.
\end{equation}
From Theorem \ref{maintheorem2} follows that, for any positive integers $n$ and $k$  such that $n-k\geq 1$ and $k>2$,  we have
\[
\mathcal{EC} _{n-k}(x_2,x_3,\dots,x_{n})\geq (k+1)\operatorname{C}_{\mathrm{Lag}}(E^{2n}(1,x_2,x_3,\dots,x_{n}), \omega_{\mathrm{std}})\]
for sufficiently large $x_2$.
\end{remark}
\begin{remark}

By the monotonicity property of the Lagrangian capacity {\cite[Section 1.2]{Cieliebak2018}}, for any positive integers $n$ and $k$  such that $n-k\geq 1$, we have
\begin{equation}\label{Lagcaplowerbound}
\mathcal{EC} _{n-k}(x_2,x_3,\dots,x_{n})\geq k\operatorname{C}_{\mathrm{Lag}}(E^{2n}(1,x_2,x_3,\dots,x_{n}),\omega_{\mathrm{std}})	\end{equation}
 for any $1\leq x_2\leq x_3\leq \dots \leq x_n\leq \infty$. It is clear from (\ref{Lagcapellip}) and (\ref{SE4}) that the lower bound (\ref{Lagcaplowerbound}) is not sharp and is much weaker than $(\ref{SEh})$.

\end{remark}
We mention that using the arguments of the proof of Theorem \ref{SEh}, one also establishes that for any positive integers $n\geq 3$, there exists $0<c<\infty$ such that 
\begin{equation}\label{SEh2}
\mathcal{EB} _{n}(x_2,x_3,\dots,x_{n})\geq (n+1)\bigg(1+\frac{1}{x_2}+\cdots+\frac{1}{x_n}\bigg)^{-1}
\end{equation}
for all  $c\leq x_2\leq x_3\leq \dots\leq x_n$. However, by Buse--Hind {\cite[Theorem 2.1]{Buse-Hindellipsoidsindimensiongreaterthanfour}}, the obstruction to embeddings of large six-dimensional ellipsoids into a six-dimensional ball coming from volumes is optimal, i.e., $\mathcal{EB} _{3}(x_2,x_3)=\sqrt[3]{x_2x_3}$ for large $x_2,x_3$. In particular, this means that the lower bound in (\ref{SEh2}) is far from being sharp and is much weaker than the lower bound coming from the volume constraint. 

Theorem \ref{maintheorem2} does not cover the case $k=2$. To address this remaining case, we adopt a different strategy. The key observation is that a symplectic embedding problem can often be reformulated either as a Lagrangian embedding problem or as a classification problem for monotone Lagrangians. This perspective offers two advantages. First, it enables the use of algebraic techniques developed for monotone Lagrangians to study symplectic embedding questions. Second, it allows one to apply the splitting argument of Cieliebak--Mohnke \cite{Cieliebak2018}, which may yield sharp symplectic embedding obstructions by employing low-degree holomorphic curves with local tangency constraints. Following this, we reformulate the stabilized symplectic embedding problem for ellipsoids as a classification problem for monotone Lagrangian tori in complex projective spaces. More precisely, we offer the following to determine the asymptotic of (\ref{function-embed}) using our approach. 

\begin{theorem}\label{maintheorem1}
Suppose there exists a symplectic embedding 
\[\Phi:\big(\bar{B}^{2}(1)\times \mathbb{C}^{n-1}, \omega_{\mathrm{std}}\big) \to \big(B^{2k}(r)\times \mathbb{C}^{n-k},\omega_{\mathrm{std}}\big) \]
	for some $r<k+1$. Consider the Lagrangian torus 
\begin{equation*}
L_\Phi:=\Phi\bigg(\overbrace{S^1\bigg(\frac{r}{k+1}\bigg)\times\cdots\times S^1\bigg(\frac{r}{k+1}\bigg)}^{n-times}\bigg)\subset B^{2k}(r)\times \mathbb{C}^{n-k}\subset  \big(\mathbb{CP}^k\times \mathbb{C}^{n-k},r\omega_{\mathrm{FS}}\oplus\omega_{\mathrm{std}}\big).
\end{equation*}
The following holds.
\begin{itemize}
	\item $L_\Phi$ is monotone in both $\big(\mathbb{CP}^k\times \mathbb{C}^{n-k},r\omega_{\mathrm{FS}}\oplus\omega_{\mathrm{std}}\big)$ and  $\big(B^{2k}(r)\times \mathbb{C}^{n-k},\omega_{\mathrm{std}}\big)$.
	\item  $L_\Phi$  has the superpotential of a monotone Clifford torus in $B^{2k}(r)\times \mathbb{C}^{n-k}$. But it is not Hamiltonian isotopic in  $B^{2k}(r)\times \mathbb{C}^{n-k}$ to any Clifford torus.
	\item The superpotential of $L_\Phi$ in $\big(\mathbb{CP}^k\times \mathbb{C}^{n-k},r\omega_{\mathrm{FS}}\oplus\omega_{\mathrm{std}}\big)$ is different from the superpotentials of the Chekanov exotic torus and the monotone Clifford torus. In particular, $L_\Phi$ is an exotic torus.
\end{itemize}
\end{theorem}

\section*{Acknowledgement}
I have benefited from discussions with several mathematicians. Particularly, I am very thankful to Richard Hind, Kyler Siegel,  Kai Cieliebak, Georgios Dimitroglou Rizell, Emmanuel Opshtein, Janko Latschev, Felix Schlenk,  Alexandru Oancea, Klaus Mohnke, Milica Dukic, Amanda Hirschi, and Gorapada Bera.

I have received financial support from the ANR grant 21-CE40-0002 COSY and the Deutsche Forschungsgemeinschaft (DFG, German Research Foundation) under Germany's Excellence Strategy-The Berlin Mathematics Research Center MATH+ (EXC-2046/1, project ID: 390685689).

\section{Holomorphic planes with ends on a skinny ellipsoid}\label{dynomicsofelipsoids}
We start with a review of some fundamental computations for standard ellipsoids. Our main reference is {\cite[Section 2.1]{MR3868228}.
\begin{definition}\label{elip}
Let $0<x_1\leq x_2\leq \dots \leq x_{n}$ be positive rationally independent reals (i.e., $\frac{x_i}{x_j}\notin  \mathbb{Q}$ for $i\neq j$), the standard irrational symplectic ellipsoid is defined by
\[E^{2n}(x_1,x_2,\dots,x_n):=\bigg\{(z_1,\dots, z_{n})\in \mathbb{C}^{n}: \sum_{i=1}^{n}\frac{  \pi|z_i|^2}{x_i}\leq 1\bigg \}.\] 
It is equipped with the standard symplectic form $\omega_{\mathrm{std}}:=\sum_{i=1}^{n}dq_i\wedge dp_i$, where $z_i=q_i+ip_i$.
\end{definition}
The standard contact form $\lambda_{\mathrm{std}}$  on the boundary $\partial E^{2n}(x_1,x_2,\dots, x_{n})$  is written as 
\[\lambda_{\mathrm{std}}:=\frac{1}{2}\sum_{i=1}^{n}(q_idp_i-p_idq_i).\]

The Reeb vector field of $\lambda_{\mathrm{std}}$ on $\partial E^{2n}(x_1,x_2,\dots, x_{n})$ is given by 
\[R_{\lambda_{\mathrm{std}}}:=2\pi\sum_{i=1}^{n}\frac{1}{x_i}\big( q_i \frac{\partial}{\partial p_i}-p_i \frac{\partial}{\partial q_i}\big).\]
The Reeb flow has precisely $n$ simple periodic orbits on  $\partial E^{2n}(x_1,x_2,\dots, x_{n})$, and they are
\[\beta_j(t):=\sqrt{\frac{x_j}{\pi}}\exp\bigg(\frac{2\pi i t}{ x_j}\bigg)e_j,\]
where $j=1,2,\dots, n$ and $e_j$ is the $j$-th vector in the canonical basis of $\mathbb{C}^n$ as a complex vector space.
We denote by  $\beta^k_j$  the $k$-fold cover of $\beta_j$. We will call  ``short orbits'' to periodic orbits of the form $\beta^k:=\beta^k_1$ for $k\in \mathbb{Z}_{\geq 1}$. 

All the orbits described are contractible. Choose a spanning disk $u: D^2\to \partial E^{2n}(x_1,\dots,x_n)$ of $\beta^{k}_i$. Every symplectic trivialization of $(u^*\xi, d\lambda_{\mathrm{std}})$ restricts to a trivialization on the boundary $((\beta^{k}_{i})^*\xi, d\lambda_{\mathrm{std}})$. There is a unique trivialization of this type up to homotopy. Let $\tau_{\mathrm{ext}}$ be one such trivialization.
\begin{theorem}[Gutt--Hutchings {\cite[Section 2.1]{MR3868228}}]\label{CZIimp}
The orbits $\beta^{k}_i$ are nondegenerate, and their Conley--Zehnder index with respect to the trivialization $\tau_{\mathrm{ext}}$ is given by 
\[\operatorname{CZ}^{\tau_{\mathrm{ext}}}(\beta^{k}_i)=n-1+2\sum_{j=1}^{n} \bigg\lfloor\frac{kx_i}{x_j}\bigg \rfloor.\]
In particular, for the $k$-fold short Reeb orbit $\beta^{k}$ and $x_2>k$ we have
\begin{equation}\label{CZI}
\operatorname{CZ}^{\tau_{\mathrm{ext}}}(\beta^k)=n-1+2k.
\end{equation}
\end{theorem}

For $\vec{x}=(1,x_2,\dots,x_n)\in \{1\}\times \mathbb{R}^{n-1}_{\geq 0}$ such that $1\leq x_2\leq \dots\leq x_n$, we set $E^{2n}(\vec{x}):=E^{2n}(1, x_2,\dots, x_n)$. Consider a closed symplectic manifold $(X,\omega)$. For any $\vec{x}$, we can find $\epsilon>0$  and a symplectic embedding $i: E^{2n}(\epsilon \vec{{x}})\to X$. The symplectic manifold $X\setminus i(E^{2n}(\epsilon \vec{{x}}))$ is a symplectic cobordism with negative boundary $(\partial E^{2n}(\epsilon \vec{{x}}), \lambda_{\mathrm{std}})$ and empty positive boundary. The symplectic completion  of $X\setminus i(E^{2n}(\epsilon \vec{{x}}))$ is obtained by gluing the negative cylindrical end  $((-\infty,0]\times \partial E^{2n}(\epsilon \vec{{x}}), d(e^r \lambda_{\mathrm{std}}))$ to $X\setminus i(E^{2n}(\epsilon \vec{{x}}))$ along $\partial E^{2n}(\epsilon \vec{{x}})$. We denote it by $\widehat{X\setminus i(E^{2n}}(\epsilon \vec{{x}}))$. The symplectic form $\omega$ extends to a symplectic form on the completion as 
	\begin{equation*}
		\widehat{\omega}:=
		\begin{cases}
\omega & \text{on } X\setminus i(E^{2n}(\epsilon \vec{{x}}))\\
			d(e^r\lambda_{\mathrm{std}}) & \text{on } (-\infty,0]\times \partial E^{2n}(\epsilon \vec{{x}}).
		\end{cases}
	\end{equation*}
	
An almost complex structure $J$ on the symplectic completion $\widehat{X\setminus i(E^{2n}}(\epsilon \vec{{x}}))$ is called $\mathrm{SFT}$-admissible if it is compatible with the symplectic form $\widehat{\omega}$ and if $J$ is $r$-translation invariant in a neighborhood of the cylindrical end, preserves $\xi$ and maps $\partial_r$ to the Reeb vector field $R_{\lambda_{\mathrm{std}}}$. We define  $\mathcal{J}(\widehat{X\setminus E^{2n}}, \widehat{\omega})$ to be the set of all SFT--admissible almost complex structure on the completion $\widehat{X\setminus i(E^{2n}}(\epsilon \vec{{x}}))$.  

Consider $J\in \mathcal{J}(\widehat{X\setminus E^{2n}}, \widehat{\omega})$, $m\in \mathbb{Z}_{\geq 1}$, and  $A\in H_2(X, i(E^{2n}(\epsilon \vec{{x}})), \mathbb{Z})= H_2(X, \mathbb{Z})$. Define  
\begin{equation*}
	\mathcal{M}^{J,s}_{\widehat{X\setminus i(E^{2n}}(\epsilon \vec{{x}})), A}(\beta^m):=\left\{
	\begin{array}{l}
		u:(\mathbb{C},i) \to (\widehat{X\setminus i(E^{2n}}(\epsilon \vec{{x}})),J),\\
		du\circ i=J\circ du ,\\
		u \text{ is asymptotic to  $\beta^m$ at $\infty$,} \\
		u \text{ represents the class $A$,}\\
		u \text{ is somewhere injective.}
	\end{array}
	\right\}\Bigg/\operatorname{Aut}(\mathbb{C},i).
\end{equation*}
By standard arguments \cite{Wendl:2016aa}, for generic $J\in \mathcal{J}(\widehat{X\setminus E^{2n}}, \widehat{\omega})$, this moduli space an oriented smooth manifold. By Theorem \ref{CZIimp}, its dimension is equal to
\[2c_1(A)-2-2m.\]
In particular, this moduli space is rigid for $m=c_1(A)-1$.

\begin{theorem}[{\cite[Section 3]{MR4332489}}]\label{keytheorem1}
Let $(X,\omega)$ be a semipositive symplectic manifold of dimension $2n$. For $m=c_1(A)-1$, the signed count of elements in  $\mathcal{M}^{J,s}_{\widehat{X\setminus i(E^{2n}}(\epsilon \vec{{x}})), A}(\beta^m)$, denoted by $\#\mathcal{M}^{J,s}_{\widehat{X\setminus i(E^{2n}}(\epsilon \vec{{x}})), A}(\beta^m)$, does not depend on the generic $J\in \mathcal{J}(\widehat{X\setminus E^{2n}}, \widehat{\omega})$, the embedding $i$, the scaling $\epsilon>0$, and the parameters $0<x_2< x_3<\dots< x_n$ defining the vector $\vec{x}$ provided that  $x_2>m$.
\end{theorem}
We have the following.
\begin{theorem}	\label{importantcount}
Let $n, k\in \mathbb{Z}_{\geq 1} $ such that $n-k\geq 1$. For $(X,\omega)=(\mathbb{CP}^k\times T^{2(n-k)},\omega_{\mathrm{FS}}\oplus\omega_{\mathrm{std}})$ and $A=[\mathbb{CP}^1]\times\{*\}]$, we have 
\[\big|\#\mathcal{M}^{J,s}_{\widehat{X\setminus i(E^{2n}}(\epsilon \vec{{x}})), A}(\beta^k)\big|= (k-1)!\]
for $x_2>k$, where $\vec{x}=(1,x_2,\dots,x_n)$.
\end{theorem}
A proof of Theorem \ref{importantcount} is given in Section \ref{proofofimportantcount}.

Consider a closed symplectic manifold $(X, \omega)$ and a point $p\in X$. Let $O(p)$ denote a small unspecified neighborhood of $p$ in $X$. Let $D\subset O(p)$ denote a local real codimension-$2$ submanifold containing the point $p$. We define 
\begin{equation*}
	\mathcal{J}_D(X,\omega):=\left\{J: 
	\begin{array}{l}
		\text{$J$ is a smooth  almost complex structure on $X$},\\
		J \text{ is compatible with $\omega$},\\
		\text{$J$ is integrable on $O(p)$} ,\\
		\text{$D$ is $J$-holomorphic.} \\
	\end{array}
	\right\}.
\end{equation*}

\begin{definition}\label{tangdef}
Let $z_0\in \mathbb{CP}^1$, $J \in \mathcal{J}_{D}(X,\omega)$ and  $u:(\mathbb{CP}^1,i)\mapsto (X,J)$ be a $J$-holomorphic sphere  with $u(z_0)=p$. Choose a holomorphic function $g:O(p)\to \mathbb{C}$ such that $g(p)=0$, $dg(p)\neq 0$,  and $D=g^{-1}(0)$. Choose a holomorphic coordinate chart $h:\mathbb{C}\to \mathbb{CP}^1$ such that $h(0)=z_0$. For $k\in \mathbb{Z}_{\geq 1}$, we say  the curve $u$ satisfies the tangency constraint $\ll \mathcal{T}_D^{k-1}p\gg$ at $z_0$  if the function $g\circ u\circ h: \mathbb{C}\to \mathbb{C}$ satisfies
\[\frac{d^i}{d^iz}(g\circ u\circ h)\bigg|_{z=0}=0, \]
for all  $i=0,1,\dots,k-1.$
\end{definition}
 This notion of tangency does not depend on the choice of the functions $h$ and $g$. It only depends on the germ of $D$ near $p$; see {\cite[Section $6$]{Cieliebak_2007}} for a proof. 
 
 By {\cite[Lemma 7.1]{Cieliebak_2007}},  the tangency constraint $\ll \mathcal{T}_D^{k-1}p\gg$ can be interpreted as a local intersection number of the image of $u$ with the divisor $D$ as follows: Choose\footnote{The symplectic form $\omega$ is exact on $O(p)$, so $u$ cannot be contained in the divisor $D$ by Stokes' theorem.} a small ball $B$ around $z_0$ such that $ u^{-1}(D)\cap B=\{z_0\}$. Smoothly perturb $u|_{B}$ away from $\partial B$ to make it transverse to $D$. The signed count of transverse intersections of $u|_{B}$ with $D$ equals $k$.

We define  $\mathcal{J}(\widehat{E^{2n}}(\vec{{x}}), \widehat{\omega}_{\mathrm{std}})$ to be the set of all SFT--admissible almost complex structure on the completion $\widehat{E^{2n}}(\vec{{x}})$. Let $D_1 \subset E^{2n}(\vec{{x}})$ denote the symplectic submanifold given by $z_1= 0$, and let $J \in \mathcal{J}(\widehat{E^{2n}}(\vec{{x}}), \widehat{\omega}_{\mathrm{std}})$ such that it agrees with the standard complex structure $J_{\mathrm{std}}$ near the origin $0\in \mathbb{C}^n$. Define
\begin{equation*}
\mathcal{M}^J_{\widehat{E^{2n}}(\vec{x})}(\beta^m)\ll \mathcal{T}_{D_1}^{m-1}0\gg:=\left\{
	\begin{array}{l}
		u:(\mathbb{C},i) \to (\widehat{E^{2n}}(\vec{{x}}),J),\\
		du\circ i=J\circ du ,\\
		u \text{ is asymptotic to  $\beta^m$ at $\infty$,} \\
		u(0)=0 \text{ and satisfies $\ll \mathcal{T}_{D_1}^{k-1}0\gg$.}\\
	\end{array}
	\right\}\Bigg/\operatorname{Aut}(\mathbb{C},0).
\end{equation*}

\begin{theorem}[{\cite[Section 4.6]{Pereira:2025aa}}, cf. {\cite[Lemma 4.1.3, Lemma 3.1.11]{MR4332489}}, {\cite[Theorem 5.2.1]{Mikhalkin:2023aa}}]\label{plan}
Let $\vec{x}=(1,x_2,\dots,x_n)\in \{1\}\times \mathbb{R}^{n-1}_{\geq 0}$ such that $1\leq x_2\leq \dots\leq x_n$ and $\epsilon>0$. Let $J \in \mathcal{J}(\widehat{E^{2n}}(\epsilon\vec{{x}}), \widehat{\omega}_{\mathrm{std}})$ be an almost complex structure such that it agrees with the standard complex structure $J_{\mathrm{std}}$ near the origin $0\in \mathbb{C}^n$ and the submanifold
\[D_i:=\big\{(z_1,z_2,\dots, z_n)\in E^{2n}(\epsilon \vec{{x}}): z_i=0\big\}\]
is $J$-holomorphic for every $i=1,2,\dots,n$. For any positive integer $m$ such that $x_2>m$, the following holds.
\begin{itemize}
	\item Every $J$-holomorphic plane $u$ in $\widehat{E^{2n}}(\epsilon \vec{x})$ that satisfies the constraint $\ll \mathcal{T}_{D_1}^{m-1}0\gg$ and is asymptotic to  a short orbit $\beta^l$ with $l \leq m$ is an $m$-fold cover of the unique embedded $J$-holomorphic plane described by $z_2=z_3=\dots=z_n=0$. In particular, $l=m$, i.e., $u$ belongs to $\mathcal{M}^J_{\widehat{E^{2n}}(\epsilon\vec{x})}(\beta^m)\ll \mathcal{T}_{D_1}^{m-1}0\gg$.
	\item The moduli space $\mathcal{M}^J_{\widehat{E^{2n}}(\epsilon\vec{x})}(\beta^m)\ll \mathcal{T}_{D_1}^{m-1}0\gg$ consists of a unique transversely cut out plane that is an $m$-fold cover of the unique embedded $J$-holomorphic plane described by $z_2=z_3=\dots=z_n=0$.
\end{itemize}
\end{theorem}
\subsection{Proof of Theorem \ref{importantcount}}\label{proofofimportantcount}
Let $n$ and $k$  be positive integers such that $n-k\geq 1$. Let $\vec{x}=(1,x_2,\dots,x_n)$, where $1\leq x_2\leq \dots \leq x_n$ and $x_2>k$. Choose $\epsilon>0$ such that $E^{2n}(\epsilon \vec{x}):=E^{2n}(\epsilon,\epsilon x_2,\dots,\epsilon x_n)\xrightarrow[]{s} \mathbb{CP}^k\times T^{2(n-k)}$. By Theorem \ref{keytheorem1}, we can assume that $E^{2n}(\epsilon \vec{x})$ lies in the complement of the symplectic hypersurface $\mathbb{CP}^{k-1}\times T^{2(n-k)}$ in $\mathbb{CP}^{k}\times T^{2(n-k)}$, i.e.,
\[(E^{2n}(\epsilon \vec{x}),\omega_{\mathrm{std}})\subset (B^{2k}(1)\times T^{2(n-k)},\omega_{\mathrm{std}}\oplus\omega_{\mathrm{std}})\subset (\mathbb{CP}^{k}\times T^{2(n-k)}, \omega_{\mathrm{FS}}\oplus \omega_{\mathrm{std}}).\]
We denote by $\widehat{W}_{\mathrm{ellip}}$ the symplectic completion of $\mathbb{CP}^{k}\times T^{2(n-k)}\setminus E^{2n}(\epsilon \vec{x})$ and set $\Omega:=\omega_{\mathrm{FS}}\oplus \omega_{\mathrm{std}}$.
\begin{itemize}
\item[\textbf{Step 01}]  Choose a generic family of $\Omega$-compatible almost complex structures $J_j$ on $\mathbb{CP}^k\times T^{2(n-k)}$ obtained via neck-stretching along the contact type hypersurface $\partial E^{2n}(\epsilon \vec{x})$. We assume that 
\begin{itemize}
\item[(a)] $J_j$ restricted to a small neighborhood of the hypersurface $\mathbb{CP}^{k-1}\times T^{2(n-k)}$ at infinity is the standard complex structure $J_{\mathrm{std}}:=J_{\mathrm{std}}\oplus J_{\mathrm{std}}$ so that hypersurface $\mathbb{CP}^{k-1}\times T^{2(n-k)}$ is $J_j$-holomorphic for all $j$;
\item[(b)] $J
_j$ restricted to the embedded ellipsoid $E^{2n}(\epsilon \vec{x})$ is such that  the submanifold
\[D_i:=\{(z_1,z_2,\dots, z_{n})\in \mathbb{C}^{n}: z_i=0\}\]
is $J_j$-holomorphic for all $j$ and each $i=1,2,\dots,n;$ and 
\item[(c)] near the center, denoted by $0$, of the ellipsoid $E^{2n}(\epsilon \vec{x})$, $J_j$ agrees with the standard complex structure $J_{\mathrm{std}}$ for all $j$. 
\end{itemize}
Let $J_{\infty}$ and $J_{\mathrm{bot}}$ denote the almost complex structures on $\widehat{W}_{\mathrm{ellip}}$ and $\widehat{E^{2n}}(\epsilon \vec{x})$, respectively, obtained as the limit of $J_j$. By construction, the hypersurface $\mathbb{CP}^{k-1}\times T^{2(n-k)}$ in $\widehat{W}_{\mathrm{ellip}}$ is $J_\infty$-holomorphic.

\item[\textbf{Step 02}] By {\cite[Corollary 3.10]{Faisal:2024}}, for each  $J_j$, there exists a $J_j$-holomorphic sphere in $\mathbb{CP}^{k}\times T^{2(n-k)}$ in the homology class $[\mathbb{CP}^1\times\{*\}]$ satisfying the constraint $\ll \mathcal{T}_{D_1}^{k-1}0\gg$, where the local symplectic divisor $D_1$ is defined by the equation $z_1= 0$ in the embedded ellipsoid $E^{2n}(\epsilon \vec{x})$.	
\item[\textbf{Step 03}] As $j\to \infty$, the sequence of $J_j$-holomorphic spheres in \textbf{Step 02} breaks into a holomorphic building $\mathbb{H}$ with its top level $\widehat{W}_{\mathrm{ellip}}$,  bottom level in $\widehat{E^{2n}}(\epsilon \vec{x})$, and some symplectization levels in $\mathbb{R}\times \partial E^{2n}(\epsilon \vec{x})$. The building $\mathbb{H}$ is of genus zero. It has index zero, i.e., the sum of the indices of its curve components equals zero. Moreover, $\mathbb{H}$ represents the homology class $[\mathbb{CP}^1\times\{*\}]$.  Since only contributions to the homology class $[\mathbb{CP}^1\times\{*\}]$ come from the top level, the top level cannot be empty. Also, the constraint $\ll \mathcal{T}_{D_1}^{k-1}0\gg$ is inherited by the portion of the building in the bottom level.

  We prove that the top level of the building in $\widehat{W}_{\mathrm{ellip}}$ consists of a somewhere injective negatively asymptotically cylindrical rigid $J_{\infty}$-holomorphic plane, denoted by $u_{\infty}$, that computes the moduli space  
\begin{equation}\label{keymoduli}
\mathcal{M}^{J_\infty}_{\widehat{W}_{\mathrm{ellip}}, [\mathbb{CP}^1\times\{*\}]}(\beta^{k}):=\left\{
	\begin{array}{l}
		u:(\mathbb{C},i) \to (\widehat{W}_{\mathrm{ellip}},J_\infty),\\
		du\circ i=J_\infty \circ du ,\\
		u \text{ is asymptotic to  $\beta^k$ at $\infty$,} \\
		u \text{ represents the class $[\mathbb{CP}^1\times\{*\}]$}
	\end{array}
	\right\}\Bigg/\operatorname{Aut}(\mathbb{C},i).
\end{equation}
Note that every curve in the moduli space (\ref{keymoduli}) is simple. The bottom level consists of a smooth positively asymptotically cylindrical $J_{\mathrm{bot}}$-plane, denoted by $u_{\mathrm{bot}}$, in $\widehat{E^{2n}}(\epsilon \vec{x})$ that inherits the constraint $\ll \mathcal{T}_{D_1}^{k-1}0\gg$. More precisely, $u_{\mathrm{bot}}$ computes the moduli space
\begin{equation}\label{keymoduli1}
\mathcal{M}^{J_{\mathrm{bot}}}_{\widehat{E^{2n}}(\epsilon \vec{x})}(\beta^{k})\ll \mathcal{T}_{D_1}^{k-1}0\gg:=\left\{
	\begin{array}{l}
		u:(\mathbb{C},i) \to \big(\widehat{E^{2n}}(\epsilon \vec{x}),J_{\mathrm{bot}}\big),\\
		du\circ i=J_{\mathrm{bot}}\circ du ,\\
		u \text{ is asymptotic to  $\beta^k$ at $\infty$,} \\
		u(0)=0 \text{ and satisfies $\ll \mathcal{T}_{D_1}^{k-1}0\gg$}\\
	\end{array}
	\right\}\Bigg/\operatorname{Aut}(\mathbb{C},0).
\end{equation}

 Moreover, there are no symplectization levels. 
	  For an illustration, see Figure \ref{Ellip3}.
	  
\begin{figure}[h]
\centering
\includegraphics[width=11cm]{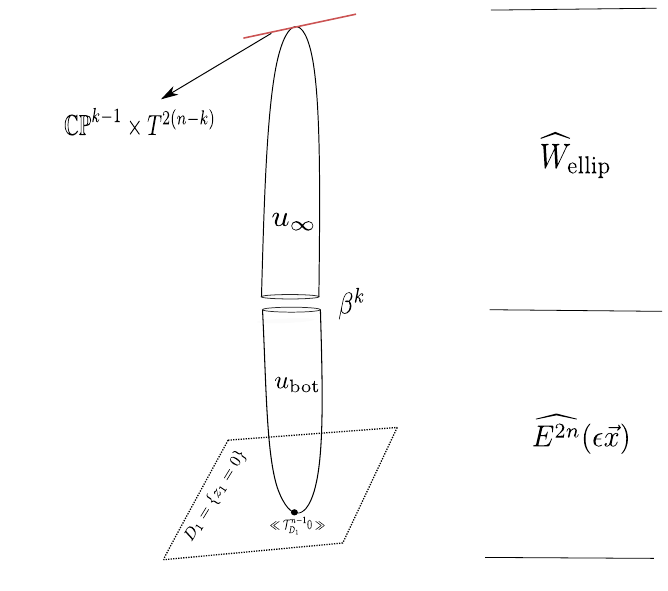}
\caption{The holomorphic building $\mathbb{H}$}\label{Ellip3}
\end{figure}

\item [\textbf{Step 04}] We start by proving that the holomorphic building $\mathbb{H}$ has no node between its non-constant curve components. In particular, the building does not have a node at the constrained marked point $0\in E^{2n}(\epsilon \vec{x})$. Since $\mathbb{H}$ has genus zero, any node between non-constant components decomposes the building $\mathbb{H}$ into two pieces $A_1$ and $A_2$. These are represented by non-constant holomorphic curves (lying possibly at different levels). By {\cite[Lemma 2.6]{Cieliebak_2005}}, we have $\Omega(A_1)>0, \Omega(A_2)>0$, where $\Omega=\omega_{\mathrm{FS}}\oplus \omega_{\mathrm{std}}$. The building $\mathbb{H}$ has exactly one intersection with the $J_\infty$-holomorphic hypersurface $\mathbb{CP}^{k-1}\times T^{2(n-k)}$  because of the positivity of intersection and the fact that $\mathbb{H}$ is the result of breaking of a holomorphic sphere that has intersection number $+1$ with $\mathbb{CP}^{k-1}\times T^{2(n-k)}$. Therefore, one of these spheres, say $A_1$, lies in the complement of $\mathbb{CP}^{k-1}\times T^{2(n-k)}$. But since $\Omega=\omega_{\mathrm{FS}}\oplus \omega_{\mathrm{std}}$ is exact in the complement of the hypersuface $\mathbb{CP}^{k-1}\times T^{2(n-k)}$, we must have $\Omega(A_1)=0$ by Stokes' theorem. This is a contradiction. We conclude that a smooth component, denoted by  $u_{\mathrm{bot}}$, in the bottom level $\widehat{E^{2n}}(\epsilon \vec{x})$ inherits the tangency constraint $\ll\mathcal{T}_{D_1}^{k-1}0\gg$.
\item [\textbf{Step 05}] Let $u_1, u_2,\dots, u_l, u_\infty$ denote the smooth connected components of the building $\mathbb{H}$ in the top level. By construction of $\mathbb{H}$ and positivity of intersection, exactly one curve component, say $u_\infty$, intersects the complex hypersurface $\mathbb{CP}^{k-1}\times T^{2(n-k)}$, and the intersection number is $+1$. Because the intersection number of $u_\infty$ with $\mathbb{CP}^{k-1}\times T^{2(n-k)}$ is $+1$, $u_\infty$ somewhere injective. The other components $u_1, u_2,\dots, u_l$ must be contained in the complement of $\mathbb{CP}^{k-1}\times T^{2(n-k)}$. 

Consider the $2$-form $\tilde{\Omega}$ on the completion of $\mathbb{CP}^k\times T^{2(n-k)}\setminus E^{2n}(\epsilon \vec{x})$,  defined by 
	\begin{equation*}
		\tilde{\Omega}:=
		\begin{cases}
\omega_{\mathrm{FS}}\oplus \omega_{\mathrm{std}} & \text{on } \mathbb{CP}^k\times T^{2(n-k)}\setminus E^{2n}(\epsilon \vec{x})\\
			d\lambda_{\mathrm{std}} & \text{on } (-\infty,0]\times \partial E^{2n}(\epsilon \vec{x}).
		\end{cases}
	\end{equation*}
The symplectic form $\omega_{\mathrm{FS}}$ is exact on the complement of the hypersurface $\mathbb{CP}^{k-1}$ in $\mathbb{CP}^k$ and $\pi_2(T^{2(n-k)})=0$. Therefore, none of $u_1, u_2,\dots, u_l$ is a closed $J_\infty$-holomorphic sphere, i.e., each of $u_1,\dots, u_l$ is a negatively asymptotically cylindrical curve with at least one negative end. It hence must have a negative $\tilde{\Omega}$-area. This is a contradiction. The conclusion is that the top level consists of exactly one smooth somewhere injective curve component, which is a punctured sphere denoted by $u_\infty$, that intersects the hypersurface $\mathbb{CP}^{k-1}\times T^{2(n-k)}$.
\item [\textbf{Step 06}] We prove that all negative punctures of $u_\infty$ are asymptotic to covers (possibly multiple) of the short Reeb orbit $\beta_1$. Suppose  $u_\infty$ has negative ends on the Reeb orbits $\beta^{m_1}_{i_1}, \beta^{m_2}_{i_2},\dots,\beta^{m_l}_{i_l}$ and assume at least one, say $\beta^{m_1}_{i_1}$, is a long orbit. The Fredholm index of $u_\infty$ in the trivialization $\tau_{\mathrm{ext}}$ (cf. Theorem \ref{CZIimp}) is
\[\operatorname{ind}(u_\infty)=(n-3)(2-l)+2c _1([\mathbb{CP}^1])-\sum_{j=1}^{l}\operatorname{CZ}^{\tau_{\mathrm{ext}}}(\beta^{m_j}_{i_j}).\] 
One can see that
\begin{equation}\label{local1}
\operatorname{ind}(u_\infty)\leq (n-3)+2c _1([\mathbb{CP}^1])-\operatorname{CZ}^{\tau_{\mathrm{ext}}}(\beta^{m_1}_{i_1}).
\end{equation}
By Theorem \ref{CZIimp}, for the long orbit $\beta^{m_1}_{i_1}$ we have
\begin{equation}\label{local2}
\operatorname{CZ}^{\tau_{\mathrm{ext}}}(\beta^{m_1}_{i_1})\geq n-1+2(\lfloor x_2\rfloor+1).
\end{equation}

Combining Equations (\ref{local1}) and (\ref{local2}) yields
	\[\operatorname{ind}(u_\infty)\leq 2(c _1([\mathbb{CP}^1\times\{*\}])-\lfloor x_2\rfloor-2).\] 
	Note that  $c _1([\mathbb{CP}^1\times\{*\}])=k+1$ and by our assumption  $x_2> k$, we have
	\[\operatorname{ind}(u_\infty)\leq 2(k-1-\lfloor x_2\rfloor) \leq -2. \]
The curve $u_\infty$ is simple and  $J_\infty$-holomorphic. We can perturb $J_\infty$ near the hypersurface $\mathbb{CP}^{k-1}\times T^{2(n-k)}$ to assume $u_\infty$ is regular. So we must have $\operatorname{ind}(u_\infty)\geq 0$. This contradicts the above estimate. Thus, all the ends of  $u_\infty$ are on short Reeb orbits.

\item[\textbf{Step 07}]  We now show that $u_\infty$ has a single negative puncture, which is asymptotic to the $m$-fold cover of the short Reeb orbit, denoted by $\beta^m$, for some positive integer $m\leq n$. By \textbf{Step 06}, we know that all negative punctures of $u_\infty$ are asymptotic to short  Reeb orbits. Suppose that  $u_\infty$ has negative ends on the Reeb orbits $\beta^{m_1}, \beta^{m_2},\dots,\beta^{m_l}$. The Fredholm index of $u_\infty$ in the trivialization $\tau_{\mathrm{ext}}$ is
\[\operatorname{ind}(u_\infty)=(n-3)(2-l)+2(n+1)-\sum_{j=1}^{l}\operatorname{CZ}^{\tau_{\mathrm{ext}}}(\beta^{m_j}).\] 
By Theorem \ref{CZIimp}, we have
\[\operatorname{CZ}^{\tau_{\mathrm{ext}}}(\beta^{m_i})\geq n-1+2m_i. \]
This implies
\[\operatorname{ind}(u_\infty)\leq (n-3)(2-l)+2(k+1)-l(n-1)-2\sum_{i=1}^{l}m_i.\]
If $l\geq 2$, then 
\[	\operatorname{ind}(u_\infty)\leq 2(k+1)-2(n-1)-2\sum_{i=1}^{2}m_i\]
Since $\sum_{i=1}^{2}m_i\geq 2$ and by our assumption $n\geq k+1$, we have 
\[\operatorname{ind}(u_\infty)\leq -2.\]
This is again a contradiction. So we must have $l=1$, i.e., $u_\infty$ has only one negative puncture.

Suppose the negative puncture of $u_\infty$ is asymptotic to $\beta^m$. By the same arguments as above, we have 
\[0\leq \operatorname{ind}(u_\infty)\leq 2(k-m).\]
This means $m\leq k$.

\item[\textbf{Step 08}] Recall from \textbf{Step 04} that the bottom level  $\widehat{E(a,b)}(\epsilon \vec{x})$ contains a smooth component, denoted by  $u_{\mathrm{bot}}$, that inherits the tangency constraint $\ll\mathcal{T}_{D_1}^{k-1}0\gg$. We show that $u_{\mathrm{bot}}$ has a single positive puncture that is asymptotic to a short Reeb orbit, i.e., to a cover of $\beta_1$.

The underlying graph of the building $\mathbb{H}$ is a tree since the building has genus zero. Suppose  $u_{\mathrm{bot}}$ has $m$ positive punctures, for some positive integer $m$. There are $m$ edges emanating from the vertex $u_{\mathrm{bot}}$ in the underlying graph. We order these edges from $1,2, \dots, m$. Let $C_i$ be the subtree emanating from the vertex $u_{\mathrm{bot}}$ along the $i$th edge. The trees $C_1,\dots, C_{k+1},\dots, C_{m}$ are topological planes with curve components in different levels.  Since the building has only one curve component in the top level, and that is  $u_\infty$, at most one of $C_1,\dots, C_{m}$, say $C_{m}$, contains $u_\infty$. By the maximum principle, each of $C_i$ must have some curve components in the top level.  Thus, we have at least $m$ smooth connected components in the top level. But by \textbf{Step 05}, there is only one curve component in the top level, namely $u_\infty$. Thus, we must have $m=1$,  i.e.,  $u_{\mathrm{bot}}$ has only one positive puncture.

Next, we prove that the positive puncture of $u_{\mathrm{bot}}$ is asymptotic to a short Reeb orbit. Suppose the Reeb orbit $\beta_j^l$ is the positive asymptotic of $u_{\mathrm{bot}}$. Recall that $\mathbb{H}$ is the limit of a sequence of spheres in the holomology class $[\mathbb{CP}^1\times\{*\}]$. Also, the top level consists of a single curve component $u_\infty$. Let  $\beta^m$ be the positive asymptotic of $u_\infty$. By \textbf{Step 07}, we have $m\leq k$. 
Since the total action of the Reeb orbits that appear as the negative asymptotics of curves in a given level in the building $\mathbb{H}$ decreases as one goes from the top to the bottom in the symplectization levels $\mathbb{R}\times \partial E^{2n}(\epsilon \vec{x})$, we have
\[ lx_j=\frac{1}{\epsilon}\int_{\beta_j^l}\lambda_{\mathrm{std}}\leq \frac{1}{\epsilon} \int_{\beta^m}\lambda_{\mathrm{std}}= m\leq  k.\]
By our assumption $x_n\geq  \dots \geq x_2>k$, so we must have $\beta_j^l=\beta_1^l$ for some $l \leq k$, where $\beta_1$ is the short Reeb orbit. 
\item[\textbf{Step 09}] By \textbf{Step 08}, $u_{\mathrm{bot}}$ in the bottom level  $\widehat{E^{2n}}(\epsilon \vec{x})$ is a plane that satisfies the tangency constraint $\ll\mathcal{T}_{D_1}^{k-1}0\gg$ and is asymptotic to a short Reeb orbit, i.e., to a cover of $\beta_1$.

 By Theorem \ref{plan},  $u_{\mathrm{bot}}$ belongs the transversely cut out moduli space given by
\begin{equation*}
\mathcal{M}^{J_{\mathrm{bot}}}_{\widehat{E^{2n}}(\epsilon \vec{x})}(\beta^{k})\ll \mathcal{T}_{D_1}^{k-1}0\gg:=\left\{
	\begin{array}{l}
		u:(\mathbb{C},i) \to \big(\widehat{E^{2n}}(\epsilon \vec{x}),J_{\mathrm{bot}}\big),\\
		du\circ i=J\circ du ,\\
		u \text{ is asymptotic to  $\beta^k$ at $\infty$,} \\
		u(0)=0 \text{ and satisfies $\ll \mathcal{T}_{D_1}^{k-1}0\gg$.}\\
	\end{array}
	\right\}\Bigg/\operatorname{Aut}(\mathbb{C},0).
\end{equation*}
 Moreover, this implies that the negative puncture of $u_\infty$ is asymptotic to $\beta^m$, for $m\geq k$. But by \textbf{Step 07} we have $m\leq k$. So we must have that the negative puncture of $u_\infty$ is asymptotic to $\beta^k$. Therefore, the curve $u_\infty$ computes the moduli space
\begin{equation*}
\mathcal{M}^{J_\infty}_{\widehat{W}_{\mathrm{ellip}}, [\mathbb{CP}^1\times\{*\}]}(\beta^{k}):=\left\{
	\begin{array}{l}
		u:(\mathbb{C},i) \to (\widehat{W}_{\mathrm{ellip}},J_\infty),\\
		du\circ i=J\circ du ,\\
		u \text{ is asymptotic to  $\beta^k$ at $\infty$,} \\
		u \text{ represents the class $[\mathbb{CP}^1\times\{*\}]$.}
	\end{array}
	\right\}\Bigg/\operatorname{Aut}(\mathbb{C},i).
\end{equation*}
By Theorem \ref{keytheorem1}, the signed count $\#\mathcal{M}^{J_\infty}_{\widehat{W}_{\mathrm{ellip}}, [\mathbb{CP}^1\times\{*\}]}(\beta^{k})$  does not depend on the generic SFT-admissible almost complex structure $J_\infty$ and the  embedded ellipsoid 
$E^{2n}(\epsilon \vec{x})$
as $x_2>k$. Also by Theorem \ref{plan}, we have 
\[\# \mathcal{M}^{J_{\mathrm{bot}}}_{\widehat{E^{2n}}(\epsilon \vec{x})}(\beta^{k})\ll \mathcal{T}_{D_1}^{k-1}0\gg=\pm 1.\]

On the other hand, gluing $u_\infty$ and  $u_{\mathrm{bot}}$ produces a curve that computes the curve count in {\cite[Corollary 3.10]{Faisal:2024}}. Therefore, by {\cite[Corollary 3.10]{Faisal:2024}}, we have 
\[
(k-1)!=\# \mathcal{M}^{J_{\mathrm{bot}}}_{\widehat{E^{2n}}(\epsilon \vec{x})}(\beta^{k})\ll \mathcal{T}_{D_1}^{k-1}0\gg \cdot \# \mathcal{M}^{J_\infty}_{\widehat{W}_{\mathrm{ellip}}, [\mathbb{CP}^1\times\{*\}]}(\beta^{k}).
\]
In particular,
\begin{equation}
|\# \mathcal{M}^{J_\infty}_{\widehat{W}_{\mathrm{ellip}}, [\mathbb{CP}^1\times\{*\}]}(\beta^{k})|= (k-1)!.
\end{equation}
\end{itemize}
This completes our proof.

\section{Enumerative descendants}\label{section3}
  Let $(X,\omega)$ be a compact monotone symplectic manifold, and let $d$ be a positive integer. Following {\cite{Cieliebak_2007, Cieliebak2018}}, we introduce a variant of Gromov--Witten invariants that we denote by $\langle \psi_{d-2}p\rangle_{X,d} ^\bullet$ and refer to as enumerative descendants of  $(X,\omega)$. The definition goes as follows.

Fix positive integers $d, k,$  and $N$. Choose a closed smooth oriented divisor $\Sigma$ in $X$ that is Poincar$\acute{\text{e}}$ dual to $Nc_1(X)$. Choose a point $p\in X\setminus \Sigma$ and local divisor $D\subset X\setminus \Sigma$ containing $p$.  Let $ \mathcal{J}_{D}(X,\omega, \Sigma)$ be the space of all  $\omega$-compatible almost complex structures that preserve the tangent spaces of $\Sigma$, i.e., $\Sigma$ is $J$-holomorphic. Let $J\in \mathcal{J}_{D}(X,\omega, \Sigma)$ and consider the moduli space
\begin{equation}\label{hamilto-perturb}
	\mathcal{M}^{J}_{X, Nd} ( \psi_{k-1}p, \Sigma)_d:=\left\{
	\begin{array}{l}
	(u,	0, z_1,z_2,\dots ,z_{Nd}),\\
		0, z_1,z_2,\dots ,z_{Nd}\in \mathbb{CP}^1,\\
		u:\mathbb{CP}^1\to (X,J),\\
		(du-K)^{0,1}=0, \ 	c_1([u])=d ,\\
		u(0)=p \text{  and satisfies $\ll \mathcal{T}_D^{k-1}p\gg$ at $0$}, \\
		u(z_i)\in \Sigma,\text{ for all } i=1,2,\dots ,Nd.
\end{array}
	\right\}\Bigg/\sim.
\end{equation}
The coherent
 perturbation $K$ (cf. {\cite[Section 6.3--4, and p. 274]{Cieliebak2018}}) in the holomorphic curve equation is given by $K=X_H\otimes\beta$ , where
\begin{itemize}
\item  $X_H$ is the Hamiltonian vector field of a time-independent Hamiltonian $H:M\to \mathbb{R}$ that is supported near the constrained point $p$. Moreover, $H$ is chosen so that $X_H$ is transverse to the local divisor $D$.
\item $\beta$ is a $1$-form on $\mathbb{CP}^1$ supported in an annulus around the point $0$. The annulus depends on the stable curve, i.e., on the position of the marked points $z_1,z_2,\dots,z_{Nd}$. It is chosen so that the points $z_i$ sit in the complement of the annulus on the side opposite to where $0$ sits.  
\end{itemize}
Note that every solution $u$ of $(du-  X_H \otimes \beta )^{0,1}=0$ is purely $J$-holomorphic near $0$, so the tangency constraint $\ll \mathcal{T}_D^{k-1}p\gg$ at $0$ is well-defined. Define $\langle \psi_{d-2}p\rangle_{X,d} ^\bullet$ to be the signed count 
\begin{equation}\label{gravi-invariant}
	\langle \psi_{d-2}p\rangle_{X,d} ^\bullet:=\frac{1}{(Nd)!}\#\mathcal{M}^{J}_{X, Nd} ( \psi_{d-2}p, \Sigma)_d.
\end{equation}
Note that a non-constant $J$-holomorphic sphere $u$ in $X$ with $c_1([u])=d$ has precisely $dN$ intersection points with $\Sigma$ due to the positivity of intersection. There are $(Nd)!$ ways to order these intersection points without repetition to produce the ordered marked point $z_1,z_2,\dots, z_{Nd}$. Therefore, we divide the right-hand side of (\ref{gravi-invariant}) by $(Nd)!$.

By {\cite[Theorem 1.2--3]{Cieliebak_2007}} or {\cite[Section 2.3]{Tonkonog:2018aa}}, the count defined by (\ref{gravi-invariant}) does not depend on $p, D, J$, $\Sigma$, and the  coherent
 perturbation $K$. 

By {\cite[Corollary 3.10]{Faisal:2024}}, for $(X,\omega)=(\mathbb{CP}^k\times T^{2m},\omega_{\mathrm{FS}}\oplus \mathbb{\omega}_{\mathrm{std}})$ and the indecomposable homology class $A=[\mathbb{CP}^1\times \{*\}]\in H_2(\mathbb{CP}^k\times T^{2m},\mathbb{Z})$, the count $\langle \psi_{k-1}p\rangle_{\mathbb{CP}^k\times T^{2m}, k+1} ^\bullet$ is equal to $(k-1)!$ when $K=0$, so we have the following.

\begin{theorem}[{\cite[Corollary 3.10]{Faisal:2024}}, {\cite[Theorem 1.2--3]{Cieliebak_2007}}]\label{countlast}
For every $k, m\in \mathbb{Z}_{\geq 1}$, we have
	\[\langle \psi_{k-1}p\rangle_{\mathbb{CP}^k\times T^{2m}, k+1} ^\bullet=(k-1)!.\]
\end{theorem}

\section{Gravitational descendants of cotangent bundles}

Let $(X,\lambda)$ be a non-degenerate Liouville domain.  We denoted by $\operatorname{SH^*_{S^1,+}}$ the the positive $S^1$-equivariant  symplectic cohomology of $(X,\lambda)$.  We briefly recall its definition following  \cite{MR3671507}. Let 
\[\operatorname{CF^*_{S^1,+}}(X)\]
be the $\mathbb{Q}$-vector space generated by the good\footnote{Let $\gamma$ be a closed Reeb orbit and $\bar{\gamma}$ be the underlying simple closed Reeb orbit. The closed Reeb orbit $\gamma$ is good if $\operatorname{CZ}^\tau(\gamma)$ and $\operatorname{CZ}^\tau(\bar{\gamma})$ have the same parity for some trivialization $\tau$. Recall that the parity of $\operatorname{CZ}^\tau(\cdot)$ does not depend on $\tau$. } Reeb orbits on $(\partial X, \lambda)$. Let $J$ be an SFT-admissible almost complex structure on $\widehat{X}$ such that $J|_{[0,\infty)\times \partial X}$ is the restriction of some SFT-admissible almost complex structure $J$ on the symplectization $(\mathbb{R}\times \partial X, d(e^r\lambda))$. The differential $\partial:\operatorname{CF^*_{S^1,+}}(X) \to \operatorname{CF^*_{S^1,+}}(X)$ is defined by 
\[\partial \gamma:=\sum_{\eta} \langle \partial \gamma, \eta \rangle \eta , \]
where $\langle \partial \gamma, \eta \rangle$ is the count of index $1$ punctured $J$-holomorphic spheres (without asymptotic markers) in the symplectization $\mathbb{R}\times \partial X$ with one positive puncture asymptotic to $\gamma$, a negative puncture asymptotic to $\eta$, and some additional negative punctures which are augmented by asymptotically cylindrical $J$-holomorphic planes in $X$.

 To achieve $\partial \circ  \partial=0$ and make the chain complex $\operatorname{CF^*_{S^1,+}}(X)$ independent (up to chain homotopy) on the choice of the almost complex structure $J$, one requires a suitable virtual perturbation scheme to define the curve count involved; for a proposal see {\cite{Chaidez:2024aa, Pardon-Contacthomologyandvirtualfundamentalcycles}}. When $(X,\lambda)$ is the unit codisk bundle of a closed Riemannian manifold that admits a metric of non-positive sectional curvature, then there are no closed Reeb orbits on $\partial X$ that are contractible in $X$, as in this case, there are no contractible closed geodesics. The differential $\partial$ counts pure holomorphic cylinders, which are unbranched by the Riemann--Hurwitz formula, interpolating between the input and output closed Reeb orbits.  In this case, one does not require any virtual perturbation scheme to define $\operatorname{CF^*_{S^1,+}}(X)$; in particular, this is the case for the unit codisk bundle $D^*T^n$ of the torus $T^n$ \cite{MR2475400, MR3671507} in which we are interested in this document.

We define $\operatorname{SH^*_{S^1,+}}(X)$ to be the homology of the chain complex $(\operatorname{CF^*_{S^1,+}}(X), \partial)$.

From now on, we assume that $(X,\lambda)$ is the unit codisk bundle $D^*L$ of a closed Riemannian manifold $L$ of dimension $n$ that admits a metric of non-positive sectional curvature. Up to an arbitrary high length (action) truncation, by {\cite[Lemma 2.2]{Cieliebak2018}}, we assume that every closed geodesic (closed Reeb orbit) has Morse index (Conley--Zehnder index) between $0$ and $n-1$. From now on, we assume that all the generators of the chain complex $\operatorname{CF^*_{S^1,+}}(D^*L)$ have actions smaller than a fixed number $A>0$. The grading on $\operatorname{CF^*_{S^1,+}}(D^*L)$ used in this document is the one used in \cite{Tonkonog:2018aa}. This is given by 
\begin{equation}\label{grading}
|\gamma_c|:=n-1-\mu(c),
\end{equation}
where $\mu(c)$ is the Morse index of the closed geodesic $c$ that lifts to the generator  $\gamma_c\in \operatorname{CF^*_{S^1,+}}(D^*L)$. For example, $\operatorname{CF^0_{S^1,+}}(D^*L)$ is generated by closed Reeb orbits that project to closed geodesics of Morse index $n-1$.

There is an isomorphism between the symplectic completion of $(D^*L, d\lambda_{\mathrm{can}})$ and the full cotangent bundle $(T^*L,d\lambda_{\mathrm{can}})$. Choose a point $p$ on the zero-section $L$ in  $T^*L$, a local divisor $D$ containing $p$, and a SFT-admissible almost complex structure on $(T^*L,d\lambda_{\mathrm{can}})$ that is integrable near $p$ such that $D$ is holomorphic.  For a choice of $k\geq 2$ generators  $\gamma_1, \dots, \gamma_k \in \operatorname{CF^0_{S^1,+}}(D^*L)$, define 
\begin{equation}\label{Hamil-perturb}
	\mathcal{M}^J_{T^*L}(\gamma_1, \dots, \gamma_k)\ll \psi_{k-2} p\gg :=\left\{
	\begin{array}{l}
		(u,	z_0, z_1,\dots ,z_k),\\
		z_0, z_1,\dots ,z_{k}\in \mathbb{CP}^1,\\
		u:\mathbb{CP}^1\setminus\{ z_1,\dots ,z_k\} \to (T^*L,J),\\
		(du-  X_H\otimes \beta)^{0,1}=0,\\
		u(z_0)=p \text{  and satisfies $\ll \mathcal{T}_D^{k-2}p\gg$ at $z_0$}, \\
		u \text{ is asymptotic to  $\gamma_i$ at $z_i$ for $i=1,\dots,k$.} 
	\end{array}
	\right\}\bigg/\sim.
\end{equation}
Here, the perturbation $X_H\otimes \beta$ in the holomorphic curve equation is chosen as in Section \ref{section3}. These moduli spaces are generically transversally cut out and are, moreover, rigid. Following {\cite[Section 4.4]{Tonkonog:2018aa}}, these rigid moduli spaces yield linear maps \[\langle \cdot|\cdots|\cdot \rangle: \operatorname{CF^0_{S^1,+}}(D^*L)^{\otimes k}\to \mathbb{Z}\] defined by
\[\langle \gamma_1|\gamma_2|\dots|\gamma_k \rangle:=\# 	\mathcal{M}^J_{T^*L}(\gamma_1, \dots, \gamma_k)\ll \psi_{k-2}p\gg\, \in \mathbb{Z}\] 
whenever $\gamma_i$ are generators of $\operatorname{CF^0_{S^1,+}}(D^*L)$ and extend linearly to the full complex. In fact, by {\cite[Theorem 4.2]{Tonkonog:2018aa}}, the maps $\langle \cdot|\cdots|\cdot \rangle$  belong to a $2$-family of linear maps $\{\psi_{m-1}^k\}_{m\geq 1, k\geq 2}$ that defines, for each integer $m\geq 1$, an $L_\infty$-algebra structure on the full complex $\operatorname{CF^*_{S^1,+}}(X)$, where $X$ belongs to a more general class of Liouville domains. These linear operations are called \textit{gravitational descendants} in {\cite{Tonkonog:2018aa}}.

By {\cite[Proposition 4.4]{Tonkonog:2018aa}}, for each $k\geq 2$ the maps 
 $\langle \cdot|\cdots|\cdot \rangle: \operatorname{CF^0_{S^1,+}}(D^*L)^{\otimes k}\to \mathbb{Z}$
descend to the cohomological level operations 
\begin{equation}\label{operations}
\langle \cdot|\cdots|\cdot \rangle: \operatorname{SH^0_{S^1,+}}(D^*L)^{\otimes k}\to \mathbb{Z}
\end{equation}
that are independant of the choices of  $p$, $D$, and compactly supported homotopies of $J$---the auxiliary data needed to define the moduli spaces (\ref{Hamil-perturb}). We will need the operations $\langle \cdot|\cdots|\cdot \rangle$ to count certain holomorphic buildings in terms of the ``Borman--Sheridan class'' in Section \ref{maintheorem2proof}.

\begin{theorem}[Descendants of the torus $T^n$, {\cite[Theorem 4.5]{Tonkonog:2018aa}}]\label{importcountall}
	For any  $k$ generators  $\gamma_1, \dots, \gamma_k \in \operatorname{SH^0_{S^1,+}}(D^*T^n)$,
	we have 
	\begin{equation*}
	\langle \gamma_1|\gamma_2|\dots|\gamma_k \rangle=
		\begin{cases}
			
			(k-2)! & \text{if }\, 0=\sum_{i=1}^{k}[\gamma_i]\in H_1(T^n,\mathbb{Z}),\\
			0 & \text{otherwise.} 
		\end{cases}
	\end{equation*}
\end{theorem}
\begin{remark}
Any curve contributing to  $\langle \gamma_1|\gamma_2|\dots|\gamma_{k} \rangle$ gives a null-homology of  $\sum_{i=1}^{k}\gamma_i$. Therefore, we must have $\langle \gamma_1|\gamma_2|\dots|\gamma_{k} \rangle=0$ whenever $\sum_{i=1}^{k}[\gamma_i]\neq 0$.
\end{remark}

\section{Proof of Theorem \ref{maintheorem2}}\label{maintheorem2proof}
We want to prove that, for any positive integers $n$ and $k>2$  such that $n-k\geq 1$, we have 
\begin{equation}\label{SEh}
\mathcal{EC} _{n-k}(x_2,x_3,\dots,x_{n})\geq (k+1)\bigg(1+\frac{1}{x_2}+\cdots+\frac{1}{x_n}\bigg)^{-1}.
\end{equation}
when $x_2$ is sufficiently large.
Suppose, on the contrary, that  there exist strictly increasing sequences $x_2^i\leq x_3^i \leq \dots\leq x_n^i$ such that $\lim_{i\to \infty}x_2^i=\infty$ and  symplectic embeddings  
\begin{equation}\label{symplecticembedding}
\Phi_i:\big(E^{2n}(1,x_2^i,\dots,x_n^i),\omega_{\mathrm{std}}\big)\to \big (B^{2k}(r_i)\times \mathbb{C}^{n-k},\omega_{\mathrm{std}} \big)
\end{equation}
for some 
\begin{equation}\label{assumption}
r_i<(k+1)\bigg(1+\frac{1}{x_2^i}+\cdots+\frac{1}{x_n^i}\bigg)^{-1}.
\end{equation}

We prove that this leads to a contradiction for large $x_2^i$ and any $k\geq 3$. The idea goes as follows. We write $\vec{x}^i:=(1,x_2^i,\dots,x_n^i)$, and observe that (\ref{symplecticembedding}) and  (\ref{assumption}) are equivelent to saying that there exist symplectic embeddings
\begin{equation}\label{explainidea}
\Phi_i:\big(E^{2n}(\frac{1}{r_i}\vec{x}^i),\omega_{\mathrm{std}}\big)\to \big (B^{2k}(1)\times \mathbb{C}^{n-k},\omega_{\mathrm{std}} \big)\subset (\mathbb{CP}^k\times \mathbb{C}^{n-k},\omega_{\mathrm{FS}}\oplus\omega_{\mathrm{std}})
\end{equation}
with
\begin{equation}\label{exlainidea2}
	\frac{1}{k+1} \bigg(1+\frac{1}{x_2^i}+\cdots+\frac{1}{x_n^i}\bigg)<\frac{1}{r_i}.
\end{equation}

For each $\frac{1}{t}\in (0,\frac{1}{r_i}]$, by restricting $\Phi_i$ to the smaller ellipsoid  $E^{2n}(\frac{1}{t}\vec{x}^i)$, we get an embedding 
\[\Phi_i:\big(E^{2n}(\frac{1}{t}\vec{x}^i),\omega_{\mathrm{std}}\big)\to \big (B^{2k}(1)\times \mathbb{C}^{n-k},\omega_{\mathrm{std}} \big)\subset (\mathbb{CP}^k\times \mathbb{C}^{n-k},\omega_{\mathrm{FS}}\oplus\omega_{\mathrm{std}}).\]
Denote by $\widehat{W}_{\mathrm{ellip}}$ the symplectic completion of $\mathbb{CP}^k\times \mathbb{C}^{n-k}\setminus \Phi_i(E^{2n}(\frac{1}{t}\vec{x}^i))$. Next, we pick an SFT-admissible almost complex structure $J_\infty$ on $\widehat{W}_{\mathrm{ellip}}$ and consider the moduli space 
\[\mathcal{M}^{J_\infty}_{\widehat{W}_{\mathrm{ellip}}, [\mathbb{CP}^1\times\{*\}]}(\beta^{k}):=\left\{
	\begin{array}{l}
		u:(\mathbb{C},i) \to (\widehat{W}_{\mathrm{ellip}},J_\infty),\\
		du\circ i=J_\infty \circ du ,\\
		u \text{ is asymptotic to  $\beta^k$ at $\infty$,} \\
		u \text{ represents the class $[\mathbb{CP}^1\times\{*\}]$}
	\end{array}
	\right\}\Bigg/\operatorname{Aut}(\mathbb{C},i).\]
By Theorem \ref{keytheorem1}, the signed count $\#\mathcal{M}^{J_\infty}_{\widehat{W}_{\mathrm{ellip}}, [\mathbb{CP}^1\times\{*\}]}(\beta^{k})$ is well-defined and does not depend on the choices involved provided that  $x^i_2>k$. In particular, this count does not depend on the scaling $\frac{1}{t}$ because for any two different values of $\frac{1}{t}$, the corresponding symplectic completions of $\mathbb{CP}^k\times \mathbb{C}^{n-k}\setminus \Phi_i(E^{2n}(\frac{1}{t}\vec{x}^i))$ are symplectomorphic. Moreover, for any $\frac{1}{t}\in (0,\frac{1}{r_i}]$, we must have \footnote{For any given vector $\vec{x}$, there are no obstructions to finding symplectic embeddings of $E^{2n}(\frac{1}{t}\vec{x}^i)$ for small $\frac{1}{t}$. We pick $\frac{1}{t}$ sufficiently
 small, choose an embedding of $E^{2n}(\frac{1}{t}\vec{x}^i)$ and compute the count by Theorem \ref{importantcount}.}
\[\#\mathcal{M}^{J_\infty}_{\widehat{W}_{\mathrm{ellip}}, [\mathbb{CP}^1\times\{*\}]}(\beta^{k})=(k-1)!\]
by Theorem \ref{importantcount}.
However, we prove that for $\frac{1}{t}$ satisfying 
\[\frac{1}{k+1} \bigg(1+\frac{1}{x_2^i}+\cdots+\frac{1}{x_n^i}\bigg)<\frac{1}{t}<\frac{1}{r_i}\]
we have 
\[\#\mathcal{M}^{J_\infty}_{\widehat{W}_{\mathrm{ellip}}, [\mathbb{CP}^1\times\{*\}]}(\beta^{k})=\pm 1.\]
That is, the existence of the embedding (\ref{explainidea}) satisfying (\ref{exlainidea2}) implies the existence of too few curves and hence a contradiction. To prove the latter claim, we follow the structure of arguments from {\cite[Section 6]{Faisal:2024}}.

Consider the Lagrangian torus
\[S^1\bigg(\frac{r_i}{k+1}\bigg)\times\cdots\times S^1\bigg(\frac{r_i}{k+1}\bigg)\subset E^{2n}(1,x_2^i,\dots,x_n^i).\]
The embeddings $\Phi_i$ yield a family of Lagrangian tori given by
\[L_{\Phi_i,r_i}:=\Phi_i\bigg(S^1\bigg(\frac{r_i}{k+1}\bigg)\times\cdots\times S^1\bigg(\frac{r_i}{k+1}\bigg)\bigg)\subset (B^{2k}(r_i)\times \mathbb{C}^{n-k},\omega_{\mathrm{std}}).\]

After compactifying the ball $(B^{2k}(r_i),\omega_{\mathrm{std}})$ to $(\mathbb{CP}^k,r_i\omega_{\mathrm{FS}})$, we get a family of Lagrangian tori given by
\[L_{\Phi_i,r_i}\subset (\mathbb{CP}^k\times \mathbb{C}^{n-k},r_i\omega_{\mathrm{FS}}\oplus\omega_{\mathrm{std}}).\]
Moreover, each $L_{\Phi_i,r_i}$ lies in the complement of the hypersurface $\mathbb{CP}^{k-1}\times \mathbb{C}^{n-k}$. Here we assume that $\omega_{\mathrm{FS}}$ integrates to $1$ on complex lines. 

Each $L_{\Phi_i,r_i}$ is a monotone Lagrangian torus in $(\mathbb{CP}^{k}\times \mathbb{C}^{n-k},r_i\omega_{\mathrm{FS}}\oplus\omega_{\mathrm{std}})$: note that there are $n$ Maslov index $2$ disks $u_1,\dots, u_n$ in the complement of the hypersurface $\mathbb{CP}^{k-1}\times \mathbb{C}^{n-k}$  with boundaries on $L_{\Phi_i,r_i}$ such that $\partial u_1,\dots, \partial u_n$ generate $H_1(L_{\Phi_i,r_i},\mathbb{Z})$. Moreover, each $u_j$ has a symplectic area equal to $r_i/(k+1)$. This implies $L_{\Phi_i,r_i}$ is monotone. 

Next, we analyze the Borman--Sheridan class of $L_{\Phi_i,r_i}$.
\begin{itemize}
	\item[\textbf{Step 01}] Fix $i$ and set $\Omega_i:=r_i\omega_{\mathrm{FS}}\oplus \omega_{\mathrm{std}}$. Cutting $\mathbb{C}^{n-k}$ by a sufficiently large lattice, we can see $L_{\Phi_i,r_i}$ as a monotone Lagrangian torus in the monotone symplectic manifold $(\mathbb{CP}^k\times T^{2(n-k)},\Omega_i)$ that lies in the complement of the hypersurface $\mathbb{CP}^{k-1}\times T^{2(n-k)}.$ 
	
\item[\textbf{Step 02}]	By Theorem \ref{countlast}, for generic $\Omega_i$-compatible almost complex structure $J$ on $\mathbb{CP}^k\times T^{2(n-k)}$ and generic $q\in\mathbb{CP}^k\times T^{2(n-k)}$, there exists a $J$-holomorphic sphere $u$ in the homology class $[\mathbb{CP}^1\times\{*\}]$ satisfying the constraint $\ll \mathcal{T}_D^{k-1}q\gg$ (cf. Definition \ref{tangdef}). Each such curve $u$ carries a Hamiltonian perturbation around the point $q$ as described in (\ref{hamilto-perturb}).

	\item[\textbf{Step 03}] Take a flat metric on $L_{\Phi_i,r_i}$. After scaling it, we can symplectically embed the codisk bundle of radius $2$, denoted by $D^*_2L_{\Phi_i,r_i}$, into $\mathbb{CP}^k\times T^{2(n-k)}$ in the complement of the hypersurface $\mathbb{CP}^{k-1}\times T^{2(n-k)}$. 
Perturb this metric according to {\cite[Lemma 2.2]{Cieliebak2018}} to a Riemannian metric $g$ such that, with respect to $g$, every closed geodesic $\gamma$ of length less than or equal to $c=r_i$ is noncontractible, nondegenerate (as a critical point of the energy functional) and satisfies 
	\[0\leq \operatorname{\mu}(\gamma)\leq n-1,\]
	where $\operatorname{\mu}(\gamma)$ denotes the Morse index of $\gamma$. Since $g$ can be chosen to be a small perturbation of the flat metric, we can ensure that the unit codisk bundle $D^*L_{\Phi_i,r_i}$ with respect to $g$ still symplectically embeds into $B^{2k}(r_i)\times T^{2(n-k)} = \mathbb{CP}^k\times T^{2(n-k)}\setminus \mathbb{CP}^{k-1} \times T^{2(n-k)}$.
	
	\textbf{Notations:} In the rest of the proof, we will denote by $W_{\mathrm{tor}}$ the symplectic cobordism $\mathbb{CP}^k\times T^{2(n-k)}\setminus D^*L_{\Phi_i,r_i}$ and by $\widehat{W}_{\mathrm{tor}}$ its symplectic completion. 

	\item[\textbf{Step 04}] Take a family of $\Omega_i$-compatible almost complex structures $J_j$ on $\mathbb{CP}^k\times T^{2(n-k)}$ that stretches the neck along the contact type hypersurface $S^*L_{\Phi_i,r_i}:=\partial D^*L_{\Phi_i,r_i} $. We assume that $J_j$ restricted to a small neighborhood of the hypersurface $\mathbb{CP}^{k-1}\times T^{2(n-k)}$ at infinity is the standard complex structure $J_{\mathrm{std}}\oplus J_{\mathrm{std}}$, so that the hypersurface $\mathbb{CP}^{k-1}\times T^{2(n-k)}$ is $J_j$-holomorphic for all $j$. Let $J_\infty$ be the almost complex structure on the symplectic completion $\widehat{W}_{\mathrm{tor}}$ obtained as the limit of $J_j$. Then the hypersurface $\mathbb{CP}^{k-1}\times T^{2(n-k)}$ in $\widehat{W}_{\mathrm{tor}}$ is $J_\infty$-holomorphic. We denote by  $J_{\mathrm{bot}}$ the almost complex structure on $T^*L_{\Phi_i,r_i}$ obtained as the limit of $J_j$.
	
\item[\textbf{Step 05}] Choose a point $q$ on $L_{\Phi_i,r_i}$ and a local symplectic divisor containing $q$. As $j\to \infty$, the sequence of $J_j$-holomorphic spheres in \textbf{Step 02} breaks into a holomorphic building $\mathbb{H}$ with its top level in    $\widehat{W}_{\mathrm{tor}}$---which is symplectomorphic to $\mathbb{CP}^k\times  T^{2(n-k)}\setminus L_{\Phi_i,r_i}$---bottom level in $T^*L_{\Phi_i,r_i}$, and some symplectization levels in $\mathbb{R}\times S^*L_{\Phi_i,r_i}$.

\item[\textbf{Step 06}] Following the arguments of {\cite[Section 6.3]{Faisal:2024}} for the monotone torus $L_{\Phi_i,r_i}$, we conclude that 
\begin{itemize}
	\item[(1)]  there are no symplectization levels in  $\mathbb{H}$;
	\item [(2)] the bottom level $T^* L_{\Phi_i,r_i}$ consists of a single smooth connected rigid punctured sphere $C_\mathrm{bot}$ with exactly $k+1$ positive punctures. Moreover, it inherits the tangency constraint  $\ll \mathcal{T}_D^{k-1}q\gg$ and the Hamiltonian perturbation supported around $q$;
	\item [(3)] The  top level that sits in $(\widehat{W}_{\mathrm{tor}}, J_{\infty})$ consists of $k+1$ asymptotically cylindrical somewhere injective rigid $J_\infty$-holomorphic planes $u_1,u_2,\dots, u_k, u_{\infty}$ with negative ends on $S^*L_{\Phi_i,r_i}$. Moreover, we have
\begin{equation}\label{symplecticarea}
		\int u_j^*\widehat{\Omega}_i=\frac{r_i}{k+1}.
	\end{equation}
for each $j=1, \dots, k,\infty$. By the monotonicity of $L_{\Phi_i,r_i}$, this implies that each of $u_1,\dots,u_k,u_\infty$ is of Maslov index $2$.
	\end{itemize}
   
\item[\textbf{Step 07}] 
Since the building is the limit of holomorphic spheres intersecting the $J_\infty$-holomorphic hypersurface $\mathbb{CP}^{k-1}\times T^{2(n-k)}$ with intersection number $+1$, this means at least $k$ planes, say $u_1,u_2,\dots, u_k$, lie in complement $B^{2k}(1)\times T^{2(n-k)}$ of $\mathbb{CP}^{k-1}\times T^{2(n-k)}$, since otherwise the total intersection number would be larger than $+1$. Here we use the fact that distinct holomorphic objects intersect positively. Moreover, $u_{\infty}$ has a simple intersection with $\mathbb{CP}^{k-1}\times T^{2(n-k)}$.

\item [\textbf{Step 08}]  Let $\gamma$ denote a closed Reeb orbit of action less or equal to $r_i$ on $S^*L_{\Phi_i,r_i}$. Moreover, assume it projects to a closed geodesic of Morse index $n-1$ on $L_{\Phi_i,r_i}$. Let  $\tau$ be a symplectic trivialization of $TT^*L_{\Phi,r_i}$. Define
\begin{equation*}
\mathcal{M}_{\widehat{W}_{\mathrm{tor}}}^{J_\infty}(\gamma):=\left\{
	\begin{array}{l}
		u:(\mathbb{C},i) \to (\widehat{W}_{\mathrm{tor}},J_{\infty}),\\
		du\circ i=J_{\infty}\circ du  ,\\
		u \text{ is asymptotic to $\gamma$ at $\infty$,} \\
		c_1^\tau(u)=1.
	\end{array}
	\right\}\bigg/\operatorname{Aut}(\mathbb{C},i).
\end{equation*}
The moduli space $\mathcal{M}_{\widehat{W}_{\mathrm{tor}}}^{J_\infty}(\gamma)$ consists of simple planes and has virtual dimension zero. Moreover, it is compact as it carries the minimal symplectic area. To explain this, note that any plane in this moduli space can be compactified to a Maslov index $2$ disk with boundary $L_{\Phi_i,r_i}$ which, by the monotonicity of $L_{\Phi_i,r_i}$ , must have symplectic area equal to $r_i/(k+1)$. A non-trivial holomorphic building that can appear as a result of degeneration in this moduli space contains two non-constant components in the top level (cf. {\cite[Theorem 6.10]{Faisal:2024}}), one of which can be compactified to a disk with boundary on $L_{\Phi_i,r_i}$.  Such a disk has a symplectic area of at least $r_i/(k+1)$ by the monotonicity of $L_{\Phi_i,r_i}$, leaving no symplectic area for the other component. We conclude that the signed count $\#\mathcal{M}_{\widehat{W}_{\mathrm{tor}}}^{J_\infty}(\gamma)$ is well-defined and does not depend on the choice of the generic almost complex structure.

The Borman--Sheridan class of $D^*L_{\Phi_i,r_i}$ in $\mathbb{CP}^k\times T^{2(n-k)}$ is the symplectic cohomology class defined by
\begin{equation}\label{Borman-Sheridan}
\mathcal{BS}(L_{\Phi_i,r_i}) :=\sum_{\gamma}\#\mathcal{M}_{\widehat{W}_{\mathrm{tor}}}^{J_\infty}(\gamma)\cdot \gamma \in \operatorname{SH}^0_{\mathrm{S^1,+}}(L_{\Phi_i,r_i}),
\end{equation}
where the sum is taken over $\gamma$ of degree zero (cf. Equation \ref{grading}). This class is independent of the choice of the almost complex structure $J_\infty$ because the moduli spaces appearing in its definition carry the minimal symplectic areas as explained above.

We note that the moduli spaces of the $J_\infty$-holomorphic planes $u_1,u_2,\dots,u_n, u_{\infty}$ in the top level of the building $\mathbb{H}$ in \textbf{Step 06} are computing the Borman--Sheridan class (\ref{Borman-Sheridan}).

\item [\textbf{Step 09}] From the steps above, it follows that under neck-stretching along the boundary of a Weinstein neighborhood of $L_{\Phi_i,r_i}$, the curves computing the count $\langle \psi_{k-1}p\rangle_{\mathbb{CP}^k\times T^{2m}, k+1} ^\bullet$ from Theorem \ref{countlast} descend to a two-level holomorphic building 
 \[\mathbb{H}=(u_1,u_2,\dots,u_n, u_{\infty}, C_{\mathrm{bot}})\]
with the top level consisting of $k+1$ somewhere injective rigid negatively asymptotically cylindrical $J_{\infty}$-holomorphic planes  $u_1,u_2,\dots,u_n, u_{\infty}$  in $\mathbb{CP}^k\times T^{2(n-k)}\setminus L_{\Phi,r_i}$, and a somewhere injective and rigid asymptotically cylindrical $J_{\mathrm{bot}}$-holomorphic sphere $C_{\mathrm{bot}}$ with $k+1$ positive punctures in the bottom level $(T^*L_{\Phi_i,r_i},d\lambda_{\mathrm{can}})$. The moduli spaces of the curves $u_1,u_2,\dots,u_n, u_{\infty}$ compute the class $\mathcal{BS}$ defined by (\ref{Borman-Sheridan}). The curve $C_{\mathrm{bot}}$ computes the  linear operations $\langle \cdot|\cdot|\cdots|\cdot\rangle $ defined by (\ref{operations}).  By standard gluing results,  this establishes a sign preserving bijection between the curve count $\langle \psi_{k-1}p\rangle_{\mathbb{CP}^k\times T^{2m}, k+1} ^\bullet$ and the count of holomorphic building of the type $(u_1,u_2,\dots,u_n, u_{\infty}, C_{\mathrm{bot}})$, up to the ordering of the end asymptotics. Therefore, by Theorem \ref{countlast}, we have  
\begin{equation}\label{rel}
\frac{1}{(k+1)!}	\langle \overbrace{ \mathcal{BS}|\mathcal{BS}|\dots|\mathcal{BS}}^{k+1 \text{ inputs}} \rangle=\langle \psi_{k-1}p\rangle_{\mathbb{CP}^k\times T^{2m}, k+1} ^\bullet=(k-1)!,
\end{equation}
where $\langle \cdot|\cdot|\cdots|\cdot\rangle $ are the linear operations defined in (\ref{operations}). This leads to the following conclusion. 
\begin{lemma}\label{count-imp}
There exists a generator $\gamma_\infty \in \operatorname{SH}^0_{\mathrm{S^1,+}}(D^*L_{\Phi_i,r_i})$ such that its coefficient in $\mathcal{BS}(L_{\Phi_i,r_i})$ given by the signed count of elements in the moduli space 
\begin{equation*}
\mathcal{M}_{\widehat{W}_{\mathrm{tor}}}^{J_\infty}(\gamma_\infty):=\left\{
	\begin{array}{l}
		u:(\mathbb{C},i) \to (\widehat{W}_{\mathrm{tor}},J_{\infty}),\\
		du\circ i=J_{\infty}\circ du  ,\\
		u \text{ is asymptotic to $\gamma_\infty$ at $\infty$,} \\
		c_1^\tau(u)=1,\\
		u\cdot[\mathbb{CP}^{k-1}\times T^{2(n-k)}]=+1 \\
	\end{array}
	\right\}\bigg/\operatorname{Aut}(\mathbb{C},i)
\end{equation*}
is equal to $\pm 1$.
\end{lemma}

\item [\textbf{Step 10}] 
Each symplectic embedding
\[
\Phi_i:\big(E^{2n}(1,x_2^i,\dots,x_n^i),\omega_{\mathrm{std}}\big)\to \big (B^{2k}(r_i)\times \mathbb{C}^{n-k},\omega_{\mathrm{std}} \big)
\]
gives a symplectic embedding
\[\Phi_i:\big(E^{2n}(1,x_2^i,\dots,x_n^i),\omega_{\mathrm{std}}\big)\to \big (\mathbb{CP}^{k}\times \mathbb{C}^{n-k},r_i\omega_{\mathrm{FS}}\oplus \omega_{\mathrm{std}} \big)
\]
and the image of this embedding lies in the complement of the hypersurface $\mathbb{CP}^{k-1}\times \mathbb{C}^{n-k}$ for every  $i\in \mathbb{Z}_{\geq 1}$. Moreover, we have
\[L_{\Phi_i,r_i}:=\Phi_i\bigg(S^1\bigg(\frac{r_i}{k+1}\bigg)\times\cdots\times S^1\bigg(\frac{r_i}{k+1}\bigg)\bigg)\subset \Phi_i\big(\operatorname{int}(E^{2n}(1,x_2^i,\dots,x_n^i))\big)\]
for every $i\in \mathbb{Z}_{\geq 1}$. 
This allows us to apply neck-stretching to the moduli space $\mathcal{M}_{\widehat{W}_{\mathrm{tor}}}^{J_\infty}(\gamma_\infty)$ in Lemma \ref{count-imp} in \textbf{Step 09} along the contact type hypersurface 
\[\partial E^{2n}(1,x_2^i,\dots,x_n^i)\subset \widehat{W}_{\mathrm{tor}}\approx\mathbb{CP}^k\times T^{2(n-k)}\setminus L_{\Phi_i,r_i}.\]

\item[\textbf{Step 11}] Take a sequence of SFT- admissible  almost complex structures $J_m$ on $\widehat{W}_{\mathrm{tor}}$ that stretches along the contact type hypersurface $\partial E^{2n}(1,x_2^i,\dots,x_n^i)$ in $\widehat{W}_{\mathrm{tor}}$. We assume $J_m$ restricted to a small neighborhood of the hypersurface $\mathbb{CP}^{k-1}\times T^{2(n-k)}$ is the standard complex structure $J_{\mathrm{std}}\oplus J_{\mathrm{std}}$ so that $\mathbb{CP}^{k-1}\times T^{2(n-k)}$ is $J_\infty$-holomorphic. Here $J_\infty$ is the almost complex structure on the symplectic completion of $\mathbb{CP}^k\times T^{2(n-k)}\setminus E^{2n}(1,x_2^i,\dots,x_n^i)$ obtained as the limit of $J_m$. Let $J_{\mathrm{bot}}$ denote the limiting SFT-admissible almost complex structure on the symplectic completion of $E^{2n}(1,x_2^i,\dots,x_n^i)\setminus D^*L_{\Phi_i,r_i}$.

\textbf{Notations:} Onwards, we will denote by $\widehat{W}_{\mathrm{tor}}^{\mathrm{ellip}}$ the symplectic completion of \[E^{2n}(1,x_2^i,\dots,x_n^i)\setminus D^*L_{\Phi_i,r_i}\]
and  by $\widehat{W}_{\mathrm{ellip}}$ the symplectic completion of $\mathbb{CP}^k\times T^{2(n-k)}\setminus E^{2n}(1,x_2^i,\dots,x_n^i)$.
\item[\textbf{Step 12}] Choose a sequence $u_m \in \mathcal{M}_{\widehat{W}_{\mathrm{tor}}}^{J_m}(\gamma_\infty)$. As $m\to \infty$, the $J_m$-holomorphic plane $u_m$ breaks into a holomorphic building $\mathbb{H}_\infty$ with the top level in $\widehat{W}_{\mathrm{ellip}}$,  the bottom level in $\widehat{W}_{\mathrm{tor}}^{\mathrm{ellip}}$, and some symplectization levels $\mathbb{R}\times \partial E^{2n}(1,x_2^i,\dots,x_n^i)$. We show that $\mathbb{H}_\infty$  consists of only two levels. The top level consists of a single degree one rigid $J_\infty$-holomorphic plane, denoted by $u_\infty$, asymptotic to the $k$-fold cover of the short Reeb orbit, denoted by $\beta^k$. The bottom level $\widehat{W}_{\mathrm{tor}}^{\mathrm{ellip}}$ consists of an unbranched rigid cylinder, denoted by $u_{\mathrm{cyl}}$, that is positively asymptotic to $\beta^k$ on $\partial E^{2n}(1,x_2^i,\dots,x_n^i)$ and negatively asymptotic to $\gamma_{\infty}$ on $S^*L_{\Phi_i,r_i}$. See Figure \ref{ellip4} for an illustration.

\begin{figure}[h]
\centering
\includegraphics[width=8cm]{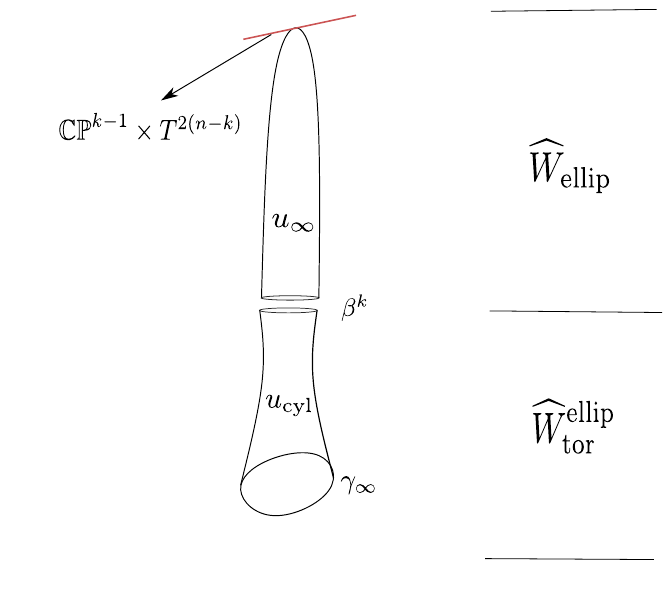}
\caption{The holomorphic building $\mathbb{H}_\infty$}\label{ellip4}
\end{figure}

\item[\textbf{Step 13}]  Let $u_1, u_2,\dots, u_l, u_\infty$ denote the smooth connected components of the building $\mathbb{H}_\infty$ in the top level. By construction, exactly one curve component, say $u_\infty$, intersects the complex hypersurface $\mathbb{CP}^{k-1}\times T^{2(n-k)}$, and the intersection number is $+1$. The other components $u_1, u_2,\dots, u_l$ are contained in the complement of $\mathbb{CP}^{k-1}\times T^{2(n-k)}$. Moreover, none of these is a closed $J_\infty$-holomorphic sphere because the symplectic form is exact on the complement of $\mathbb{CP}^{k-1}\times T^{2(n-k)}$ in $\widehat{W}_{\mathrm{ellip}}$. Because the intersection number of $u_\infty$ with $\mathbb{CP}^{k-1}\times T^{2(n-k)}$ is $+1$, $u_\infty$ somewhere injective. Moreover, the bottom level $\widehat{W}^{\mathrm{ellip}}_{\mathrm{tor}}$ contains a smooth connected curve, denoted by $u_{\mathrm{cyl}}$,  with some positive punctures asymptotic to closed Reeb orbits on $\partial E^{2n}(1,x_2^i,\dots,x_n^i)$ and a negative puncture asymptotic to the closed Reeb $\gamma_\infty$ on $S^*L_{\Phi_i,r_i}$ defined in Lemma \ref{count-imp}.

\item[\textbf{Step 14}]  The energy of $\mathbb{H}_\infty$ is given by Equation (\ref{symplecticarea}). This, together with our assumption (\ref{assumption}), implies that
\begin{equation}\label{symengerestem}
0<\int u_\infty^*\tilde{\Omega}^i_{\mathrm{ellip}}+\sum_{i=1}^{l}\int u_i^*\tilde{\Omega}^i_{\mathrm{ellip}}\leq  \frac{r_i}{k+1}<1.
\end{equation}
Here,
\begin{equation*}
	\tilde{\Omega}^i_{\mathrm{ellip}}:=
	\begin{cases}
		r_i\omega_{\mathrm{FS}}\oplus \omega_{\mathrm{std}}& \text{on } \mathbb{CP}^k\times T^{2(n-k)}\setminus E^{2n}(1,x_2^i,\dots,x_n^i),\\
		d\lambda_{\mathrm{std}}|_{\partial E^{2n}(1,x_2^i,\dots,x_n^i)} & \text{on } (-\infty,0]\times \partial E^{2n}(1,x_2^i,\dots,x_n^i).
	\end{cases}
\end{equation*}
The curve components $u_1, u_2,\dots, u_l$ are contained in the complement of $\mathbb{CP}^{k-1}\times T^{2(n-k)}$, where the 2-form $	\tilde{\Omega}^i_{\mathrm{ellip}}$ is exact, and are negatively asymptotic to closed Reeb orbits on $\partial E^{2n}(1,x_2^i,\dots,x_n^i)$.  The minimal period of a closed Reeb orbit on $\partial E^{2n}(1,x_2^i,\dots,x_n^i)$ is $1$. Thus,
\[\sum_{i=1}^{l}\int u_i^*\tilde{\Omega}^i_{\mathrm{ellip}}\geq l.\]
This contradicts (\ref{symengerestem}) if $l\geq 1$. This means $l=0$, i.e., the  curve components $u_1, u_2,\dots, u_l$ do not exist. So the top level consists of a single simple smooth connected curve $u_\infty$.
\item[\textbf{Step 15}]
We prove that all negative punctures of $u_\infty$ are asymptotic to covers (possibly multiple) of the short Reeb orbit $\beta_1$. Suppose  $u_\infty$ has negative ends on the Reeb orbits $\beta^{m_1}_{i_1}, \beta^{m_2}_{i_2},\dots,\beta^{m_l}_{i_l}$ and assume at least one, say $\beta^{m_1}_{i_1}$, is not a short orbit \footnote{An orbit that is a cover (possibly multiple) of $\beta_1$ is called a short orbit.}. The Fredholm index of $u_\infty$ in the trivialization $\tau_{\mathrm{ext}}$ (cf. Theorem \ref{CZIimp}) is
\[\operatorname{ind}(u_\infty)=(n-3)(2-l)+2c _1([\mathbb{CP}^1])-\sum_{j=1}^{l}\operatorname{CZ}^{\tau_{\mathrm{ext}}}(\beta^{m_j}_{i_j}).\] 
One can see that
\begin{equation}\label{local1}
\operatorname{ind}(u_\infty)\leq (n-3)+2c _1([\mathbb{CP}^1])-\operatorname{CZ}^{\tau_{\mathrm{ext}}}(\beta^{m_1}_{i_1})..
\end{equation}
By Theorem \ref{CZIimp}, for the long orbit $\beta^{m_1}_{i_1}$ we have
\begin{equation}\label{local2}
\operatorname{CZ}^{\tau_{\mathrm{ext}}}(\beta^{m_1}_{i_1})\geq n-1+2(\lfloor x_2^i\rfloor+1).
\end{equation}

Combining Equations (\ref{local1}) and (\ref{local2}) yields
	\[\operatorname{ind}(u_\infty)\leq 2(c _1([\mathbb{CP}^1])-\lfloor x_2^i\rfloor-2).\] 
For  $x_2^i\geq k$, using  $c _1([\mathbb{CP}^1])=k+1$, we have
	\[\operatorname{ind}(u_\infty)\leq 2(k-1-\lfloor x_2^i\rfloor) \leq -2. \]
The curve $u_\infty$ is simple and  $J_\infty$-holomorphic. We can perturb $J_\infty$ near the hypersurface $\mathbb{CP}^{k-1}\times T^{2(n-k)}$ to assume $u_\infty$ is regular. So we must have $\operatorname{ind}(u_\infty)\geq 0$. This contradicts the above estimate on its index. We conclude that all the ends of  $u_\infty$ are on short Reeb orbits.

\item[\textbf{Step 16}]  Next we prove that $u_\infty$ has a single negative puncture that is asymptotic to the $m$-fold cover of the short Reeb orbit, denoted by $\beta^m$, for some positive integer $m\leq k$. By \textbf{Step 15}, all negative punctures of $u_\infty$ are asymptotic to short  Reeb orbits. Suppose  $u_\infty$ has negative ends on the Reeb orbits $\beta^{m_1}, \beta^{m_2},\dots,\beta^{m_l}$. The Fredholm index of $u_\infty$ in the trivialization $\tau_{\mathrm{ext}}$ is
\[\operatorname{ind}(u_\infty)=(n-3)(2-l)+2(k+1)-\sum_{j=1}^{l}\operatorname{CZ}^{\tau_{\mathrm{ext}}}(\beta^{m_j}).\] 
By Theorem \ref{CZIimp}, we have
\[\operatorname{CZ}^{\tau_{\mathrm{ext}}}(\beta^{m_i})\geq n-1+2m_i. \]
This implies
\[\operatorname{ind}(u_\infty)\leq (n-3)(2-l)+2(k+1)-l(n-1)-2\sum_{i=1}^{l}m_i.\]
If $l\geq 2$, then 
\[	\operatorname{ind}(u_\infty)\leq 2(k+1)-2(n-1)-2\sum_{i=1}^{2}m_i\]
Since by our assumption $n\geq k+1$ and $\sum_{i=1}^{2}m_i\geq 2$, we have 
\[\operatorname{ind}(u_\infty)\leq -2.\]
This is again a contradiction. So we must have $l=1$, i.e., $u_\infty$ has only one negative puncture.

Suppose the negative puncture of $u_\infty$ is asymptotic to $\beta^m$. By the same arguments as above, we have 
\[0\leq \operatorname{ind}(u_\infty)\leq 2(k-m).\]
This means $m\leq k$.

\item[\textbf{Step 17}] Recall from \textbf{Step 13} that the bottom level $\widehat{W}^{\mathrm{ellip}}_{\mathrm{tor}}$ contains a smooth connected curve, denoted by $u_{\mathrm{cyl}}$,  with some positive punctures asymptotic to closed Reeb orbits on $\partial E^{2n}(1,x_2^i,\dots,x_n^i)$ and a negative puncture asymptotic to the closed Reeb $\gamma_\infty$ on $S^*L_{\Phi_i,r_i}$ defined in Lemma \ref{count-imp}. We show that $u_{\mathrm{cyl}}$ has a single positive puncture.

The underlying graph of the building $\mathbb{H}_\infty$ is a tree since the building has genus zero. Suppose  $u_{\mathrm{cyl}}$ has $m$ positive punctures, for some positive integer $m$. There are $m$ edges emanating from the vertex $u_{\mathrm{cyl}}$ in the underlying graph. We order these edges from $1,2, \dots, m$. Let $C_i$ be the subtree emanating from the vertex $u_{\mathrm{cyl}}$ along the $i$th edge. The trees $C_1,\dots, C_{k+1},\dots, C_{m}$ are topological planes with curve components in different levels.  Since the building has only one curve component in the top level, and that is  $u_\infty$, at most one of $C_1,\dots, C_{m}$, say $C_{m}$, contains $u_\infty$. By the maximum principle, each of $C_i$ must have some curve components in the top level.  Thus, we have at least $m$ smooth connected components in the top level. But by \textbf{Step 14}, there is only one curve component in the top level, namely $u_\infty$. Thus, we must have $m=1$,  i.e.,  $u_{\mathrm{cyl}}$ has only one positive puncture. We conclude that $u_{\mathrm{cyl}}$ is a cylinder with a positive and a negative end.

\item[\textbf{Step 18}] Next, we prove that the positive puncture of $u_{\mathrm{cyl}}$ is asymptotic to $\beta^k$, the $k$-fold cover of the short Reeb orbit $\beta_1$. Suppose the positive puncture of $u_{\mathrm{cyl}}$ is asymptotic to $\beta^l$, for some positive integer $l$. The cylinder $u_{\mathrm{cyl}}$ in  $\widehat{W}_{\mathrm{tor}}^{\mathrm{ellip}}\approx \widehat{E^{2n}}(1,x_2^i,\dots,x_n^i)\setminus L_{\Phi_i,r_i}$ can be compactified to a smooth half-cylinder $\bar{u}_{\mathrm{cyl}}:[0,\infty)\to \widehat{E^{2n}}(1,x_2^i,\dots,x_n^i)\setminus L_{\Phi_i,r_i}$ with $\bar{u}_{\mathrm{cyl}}(\{0\}\times S^1)\subset L_{\Phi_i,r_i}$. By construction,  the boundary of $u_{\mathrm{cyl}}$ on $L_{\Phi_i,r_i}$ bounds a disk of symplectic area $r_ik/(k+1)$, so 
\[0\leq \int \bar{u}_{\mathrm{cyl}}^*\tilde{\omega}_{\mathrm{std}}=\int_{\beta^l}\lambda_{\mathrm{std}}-\frac{r_i\, k}{k+1} \]
where $\tilde{\omega}_{\mathrm{std}}$ is the exact $2$-form
\begin{equation*}
	\tilde{\omega}_{\mathrm{std}}:=
	\begin{cases}
		d\lambda_{\mathrm{std}}|_{\partial E^{2n}(1,x_2^i,\dots,x_n^i)}& \text{on } [0,\infty)\times \partial E^{2n}(1,x_2^i,\dots,x_n^i),\\
		\omega_{\mathrm{std}} & \text{on } E^{2n}(1,x_2^i,\dots,x_n^i)\setminus L_{\Phi_i,r_i}.
	\end{cases}
\end{equation*}
Thus, we have 
\[\frac{r_i\, k}{k+1}\leq\int_{\beta^l}\lambda_{\mathrm{std}}=l. \]
Since $l$ is an integer and for sufficiently large $x_2^i$ we can choose $r_i$ 
 sufficiently closed to $k+1$ by our assumption (\ref{assumption}), therefore  $l\geq k.$

Recall that $u_\infty$ is asymptotic to $\beta^m$ for some $m\leq k$. The total action of Reeb orbits decreases as one goes from top to bottom along the building in the symplectization levels  $\mathbb{R}\times  \partial E^{2n}(1,x_2^i,\dots,x_n^i)$, so we must have $k\geq m\geq l\geq k$. Putting everything together, we obtain $m=l=k$.

The conclusion is that $u_\infty$ is negatively asymptotic to $\beta^k$ and the cylinder $u_{\mathrm{cyl}}$ is positively asymptotic to $\beta^k$. Moreover, we have $\operatorname{ind}(u_\infty)=0$ generically.

\item[\textbf{Step 19}]  The index of the building $\mathbb{H}_\infty$ is zero, so 
\[\operatorname{ind}(u_{\mathrm{cyl}})+\underbrace{\operatorname{ind}(u_\infty)}_{=0}=0.\]
This implies $\operatorname{ind}(u_{\mathrm{cyl}})=\operatorname{ind}(u_\infty)=0$. By the Riemann--Hurwitz formula, $u_{\mathrm{cyl}}$ is an unbranched cylinder.
\item[\textbf{Step 20}] From the steps above, it follows that under neck-stretching along  $\partial E^{2n}(1,x_2^i,\dots,x_n^i)$ in $\widehat{W}_{\mathrm{tor}}$, any curve computing the rigid moduli space $\mathcal{M}_{\widehat{W}_{\mathrm{tor}}}^{J_\infty}(\gamma_\infty)$ described in Lemma \ref{count-imp} in \textbf{Step 09},  descends to a two-level holomorphic building $\mathbb{H}_{\infty}=(u_\infty, u_{\mathrm{cyl}})$ with the top level consisting of a somewhere injective rigid negatively asymptotically cylindrical degree one $J_{\infty}$-holomorphic plane  $u_\infty$  in the symplectic completion of $\mathbb{CP}^k\times T^{2(n-k)}\setminus E^{2n}(1,x_2^i,\dots,x_n^i)$, denoted by $\widehat{W}_{\mathrm{ellip}}$, and an unbranched rigid asymptotically cylindrical $J_{\mathrm{bot}}$-holomorphic cylinder $u_{\mathrm{cyl}}$ in the symplectic completion of $E^{2n}(1,x_2^i,\dots,x_n^i)\setminus D^*L_{\Phi_i,r_i}$, denoted by $\widehat{W}^{\mathrm{ellip}}_{\mathrm{tor}}$. More precisely, $u_\infty$ computes the moduli space defined by
\[
\mathcal{M}^{J_\infty}_{\widehat{W}_{\mathrm{ellip}}, [\mathbb{CP}^1\times \{*\}]}(\beta^{k}):=\left\{
	\begin{array}{l}
		u:(\mathbb{C},i) \to (\widehat{W}_{\mathrm{ellip}},J_\infty),\\
		du\circ i=J_\infty \circ du ,\\
		u \text{ is asymptotic to  $\beta^k$ at $\infty$,} \\
		u \text{ represents the class $[\mathbb{CP}^1\times \{*\}]$.}
	\end{array}
	\right\}\Bigg/\operatorname{Aut}(\mathbb{C},i).
\]
By Theorem \ref{keytheorem1}, the signed count $\#\mathcal{M}^{J_\infty}_{\widehat{W}_{\mathrm{ellip}}, [\mathbb{CP}^1]}(\beta^{k})$  does not depend on the generic SFT-admissible almost complex structure $J_\infty$ and the  the embedded ellipsoid 
\[E^{2n}(1,x_2^i,\dots,x_n^i)\]
provided that  $x^i_2$ is sufficiently large.

Also the cylinder $u_{\mathrm{cyl}}$ computes the moduli space defined by
\begin{equation}\label{Hamiltonianpertur}
\mathcal{M}^{J_{\mathrm{bot}}}_{\widehat{W}^{\mathrm{ellip}}_{\mathrm{tor}}}(\beta^{k},\gamma_{\infty}):=\left\{
	\begin{array}{l}
		u:(\mathbb{R}\times S^1,i) \to (\widehat{W}^{\mathrm{ellip}}_{\mathrm{tor}},J_{\mathrm{bot}}),\\
		du\circ i=J\circ du ,\\
		u \text{ is asymptotic to  $\beta^k$ at $\infty$,} \\
		u \text{ is asymptotic to  $\gamma_\infty$ at $-\infty$,}
	\end{array}
	\right\}\Bigg/\operatorname{Aut}(\mathbb{R}\times S^1,i).
\end{equation}
The $\tilde{\omega}_{\mathrm{std}}^{\mathrm{up}}$-area of any cylinder $u$ in this moduli space is given by 
\begin{equation}\label{symplecticareaest}
0\leq \int u^*\tilde{\omega}_{\mathrm{std}}^{\mathrm{up}}=\int_{\beta^k}\lambda_{\mathrm{std}}-\frac{r_i\, k}{k+1}= k-\frac{r_i\, k}{k+1}.	
\end{equation}
where $\tilde{\omega}_{\mathrm{std}}^{\mathrm{up}}$ is the exact $2$-form
\begin{equation*}
	\tilde{\omega}_{\mathrm{std}}^{\mathrm{up}}:=
	\begin{cases}
		d\lambda_{\mathrm{std}}|_{\partial E^{2n}(1,x_2^i,\dots,x_n^i)}& \text{on } [0,\infty)\times \partial E^{2n}(1,x_2^i,\dots,x_n^i),\\
		\omega_{\mathrm{std}} & \text{on } E^{2n}(1,x_2^i,\dots,x_n^i)\setminus D^*L_{\Phi_i,r_i},\\
		d(e^r\lambda_{\mathrm{can}}) & \text{on } (-\infty,0]\times S^*L_{\Phi _i,r_i}.
	\end{cases}
\end{equation*}
Since $x_2^i$ is large, we can choose $r_i$ close to $k+1$ following our assumption (\ref{assumption}). This implies that the symplectic area given by (\ref{symplecticareaest}) is very small. Consequently, the moduli space (\ref{Hamiltonianpertur})---as well as its parametric version with respect to the almost complex structure---is compact.

 Using a suitable transversality scheme, we can assume the moduli space (\ref{Hamiltonianpertur}) is transversely cut out. For instance, following {\cite[Section 6]{Cieliebak2018}}, one can introduce a coherent Hamiltonian perturbation supported in the interior
 of 
 \[E^{2n}(1,x_2^i,\dots,x_n^i)\setminus D^*L_{\Phi_i,r_i},\] to the moduli space apearing in Lemma \ref{count-imp} in \textbf{Step 09}. Under neck-streching along $\partial E^{2n}(1,x_2^i,\dots,x_n^i)$, the moduli space (\ref{Hamiltonianpertur}) inherits this Hamiltonian perturbation. For a generic such perturbation, the moduli space (\ref{Hamiltonianpertur}) is transversely cut out by the holomorphic curve equation. The conclusion is that the integer signed count  
$\# \mathcal{M}^{J_{\mathrm{bot}}}_{\widehat{W}^{\mathrm{ellip}}_{\mathrm{tor}}}(\beta^{k},\gamma_{\infty})$ is well-defined and
does not depend on the auxiliary choices, such as the generic SFT-admissible almost complex structure $J_{\mathrm{bot}}$.

On the other hand, gluing $u_\infty$ and  $u_{\mathrm{cyl}}$ produces a curve that computes the count in the moduli space $\mathcal{M}_{\widehat{W}_{\mathrm{tor}}}^{J_\infty}(\gamma_\infty)$ described in Lemma \ref{count-imp} in \textbf{Step 09}. Therefore, by Lemma \ref{count-imp}, we have 
\[
\pm 1=\# \mathcal{M}_{\widehat{W}_{\mathrm{tor}}}^{J_\infty}(\gamma_\infty)=\#\mathcal{M}^{J_{\mathrm{bot}}}_{\widehat{W}^{\mathrm{ellip}}_{\mathrm{tor}}}(\beta^{k},\gamma_{\infty})\cdot \# \mathcal{M}^{J_\infty}_{\widehat{W}_{\mathrm{ellip}}, [\mathbb{CP}^1]}(\beta^{k}).
\]
In particular,
\begin{equation}\label{ellicountlast}
|\# \mathcal{M}^{J_\infty}_{\widehat{W}_{\mathrm{ellip}}, [\mathbb{CP}^1]}(\beta^{k})|= 1.
\end{equation}
\end{itemize}
This contradicts Theorem \ref{importantcount} for $k\geq 3$. This completes our proof.
\section{Proof of Theorem \ref{maintheorem1}}
Suppose there is an $r<k+1$ for which there exists a symplectic embedding
\begin{equation}\label{symplecticembeddinglast}
\Phi:\bar{B}^2(1)\times \mathbb{C}^{n-1}\xrightarrow[]{s} B^{2k}(r)\times \mathbb{C}^{n-k}.
\end{equation}
For $t\in [r,k+1]$, consider the Lagrangian torus
\[S^1\bigg(\frac{t}{k+1}\bigg)\times\cdots\times S^1\bigg(\frac{t}{k+1}\bigg)\subset \bar{B}^2(1)\times \mathbb{C}^{n-1}.\]
The embedding $\Phi$ yields a family of Lagrangian tori given by
\[L_{\Phi,t}:=\Phi\bigg(S^1\bigg(\frac{t}{k+1}\bigg)\times\cdots\times S^1\bigg(\frac{t}{k+1}\bigg)\bigg)\subset (B^{2k}(r)\times \mathbb{C}^{n-k},\omega_{\mathrm{std}}).\]

After compactifying the ball $(B^{2k}(r),\omega_{\mathrm{std}})$ to $(\mathbb{CP}^k,r\omega_{\mathrm{FS}})$, we get a family of Lagrangian tori given by
\[L_{\Phi,t}\subset (\mathbb{CP}^k\times \mathbb{C}^{n-k},r\omega_{\mathrm{FS}}\oplus\omega_{\mathrm{std}}).\]
Moreover, for each $t$ the torus $L_{\Phi,t}$ lies in the complement of the hypersurface $\mathbb{CP}^{k-1}\times \mathbb{C}^{n-k}.$ 

We show that $L_{\Phi,r}$ is a monotone torus that is Hamiltonian isotopic to neither the Clifford torus nor to the Chekanov--Schlenk exotic torus. Note that there are $n$ Maslov index $2$ disks $u_1,\dots, u_n$ in the complement of the hypersurface $\mathbb{CP}^{k-1}\times \mathbb{C}^{n-k}$  with boundaries on $L_{\Phi,r}$ such that $\partial u_1,\dots, \partial u_n$ generate $H_1(L_{\Phi,r},\mathbb{Z})$. Moreover, each $u_i$ has a symplectic area equal to $r/(k+1)$. This, in particular, means that $L_{\Phi,r}$ is monotone and does not belong to the Hamiltonian isotopy class of the Chekanov--Schlenk exotic torus. 

The monotone torus $L_{\Phi,r}$  has the superpotential of a monotone Clifford torus in $B^{2k}(r)\times \mathbb{C}^{n-k}$. To show this, note that for sufficiently large $S>0$ we have 
\[L_{\Phi,r}\subset \Phi(E^{2n}(1,S,\dots,S))\subset B^{2k}(r)\times \mathbb{C}^{n-k}.\]
It is enough to show that any Maslov index $2$ disk contributing to the superpotential is contained in the embedded ellipsoid $\Phi(E^{2n}(1,S,\dots,S))$. Suppose on the contrary that this is not the case, we perform neck-streching along the contact type hypersurface  $\Phi(\partial E^{2n}(1,S,\dots,S))$ to produce a negatively punctured holomorphic curve (possibly with no punctures) in the symplectic completion of $B^{2k}(r)\times \mathbb{C}^{n-k}\setminus\Phi(E^{2n}(1,S,\dots,S))$. But no such curve can exist. Because the cobordism $B^{2k}(r)\times \mathbb{C}^{n-k}\setminus \Phi(E^{2n}(1,S,\dots,S))$ is exact, so any curve with negative punctures or no punctures at all in the symplectic completion of $B^{2k}(r)\times \mathbb{C}^{n-k}\setminus \Phi(E^{2n}(1,S,\dots,S))$ will have negative or zero $\tilde{\omega}_{\mathrm{std}}$-area, respectively, where
\begin{equation*}
	\tilde{\omega}_{\mathrm{std}}:=
	\begin{cases}
		 \omega_{\mathrm{std}}& \text{on } B^{2k}(r)\times \mathbb{C}^{n-k}\setminus \Phi( E^{2n}(1,S,\dots, S)),\\
		d\lambda_{\mathrm{std}}|_{\Phi(\partial E^{2n}(1,S,\dots, S))} & \text{on } (-\infty,0]\times \Phi(\partial E^{2n}(1,S,\dots, S)).
	\end{cases}
\end{equation*}

In what follows, we prove that $L_{\Phi,r}$ is not Hamiltonian isotopic to the monotone Clifford torus in $\mathbb{CP}^k\times \mathbb{C}^{n-k}$.

\begin{itemize}
	\item[\textbf{Step 01}] 
		Set $\Omega:=r\omega_{\mathrm{FS}}\oplus \omega_{\mathrm{std}}$. By Theorem \ref{countlast}, for generic $\Omega$-compatible almost complex structure $J$ on $\mathbb{CP}^k\times T^{2(n-k)}$, there exists a $J$-holomorphic sphere in the homology class $[\mathbb{CP}^1\times\{*\}]$ passing through a generic point $p\in \mathbb{CP}^k\times T^{2(n-k)}$ and tangent of order $k-1$ to a local symplectic divisor containing $p$ (cf. Definition \ref{tangdef}). Moreover, each such curve $u$ carries a Hamiltonian perturbation around the point $p$ as described in (\ref{hamilto-perturb}).	
	
	\item[\textbf{Step 02}] 
	Take a flat metric on $L_{\Phi,t}$. After scaling it, we can symplectically embed the codisk bundle of radius $2$, denoted by $D^*_2L_{\Phi,r}$, into $\mathbb{CP}^k\times T^{2(n-k)}$ in the complement of the hypersurface $\mathbb{CP}^{k-1}\times T^{2(n-k)}$. 
Perturb this metric according to {\cite[Lemma 2.2]{Cieliebak2018}} to a Riemannian metric $g$ such that, with respect to $g$, every closed geodesic $\gamma$ of length less than or equal to $r$ is noncontractible, nondegenerate (as a critical point of the energy functional) and satisfies 
	\[0\leq \operatorname{\mu}(\gamma)\leq n-1,\]
	where $\operatorname{\mu}(\gamma)$ denotes the Morse index of $\gamma$. Since $g$ can be chosen to be a small perturbation of the flat metric, we can ensure that the unit codisk bundle $D^*L_{\Phi,r}$ with respect to $g$ still symplectically embeds into $B^{2k}(r)\times T^{2(n-k)} = \mathbb{CP}^k\times T^{2(n-k)}\setminus \mathbb{CP}^{k-1} \times T^{2(n-k)}$.
	
	\item[\textbf{Step 03}] Take a family of $\Omega$-compatible almost complex structures $J_i$ on $\mathbb{CP}^k\times T^{2(n-k)}$ that stretches along the contact hypersurface $S^*L_{\Phi,t}$ in $\mathbb{CP}^k\times T^{2(n-k)}$, where $S^*L_{\Phi,t}$ is the unit cosphere bundle of $L_{\Phi,t}$. We assume $J_i$ restricted to a small neighborhood of the hypersurface $\mathbb{CP}^{(k-1)}\times T^{2(n-k)}$ at infinity is the standard complex structure $J_{\mathrm{std}}\oplus J_{\mathrm{std}}$ so that hypersurface $\mathbb{CP}^{k-1}\times T^{2(n-k)}$ is $J_i$-holomorphic. Let $J_\infty$ be the almost complex structure on the symplectic completion of $\mathbb{CP}^k\times T^{2(n-k)}\setminus D^*L_{\Phi,t}$ obtained as the limit of $J_i$. Then the hypersurface $\mathbb{CP}^{k-1}\times T^{2(n-k)}$ is $J_\infty$-holomorphic. We denote by  $J_{\mathrm{bot}}$ the almost complex structure on $T^*L_{\Phi,r}$ obtained as the limit of $J_i$.

	\item[\textbf{Step 04}]  As $k\to \infty$, by SFT compactness theorem \cite{Cieliebak_2005, MR2026549}, the $J_k$-holomorphic sphere described in \textbf{Step 01} breaks into a holomorphic building with  top level in  $\mathbb{CP}^k\times T^{2(n-k)}\setminus L_{\Phi,t}$, bottom level in $T^*L_{\Phi,t}$, and some intermediate symplectization levels in $\mathbb{R}\times S^*L_{\Phi,t}$.
	
	\item[\textbf{Step 05}] Let $D_1, D_2,\dots, D_m$ be the smooth connected components of the building in the top level $\mathbb{CP}^k\times T^{2(n-k)}\setminus L_{\Phi,t}$. The building is the limit of holomorphic spheres of symplectic area $r$. Therefore, 
	\begin{equation}\label{embed-est}
		\sum_{j=1}^{m}\int D_j^*\Omega=r.
	\end{equation}

 There are at least $k+1$ $J_\infty$-holomorphic disks in the top level $\mathbb{CP}^k\times T^{2(n-k)}\setminus L_{\Phi,t}$, by {\cite[Lemma 6.4]{Faisal:2024}}. Denote these by $D_1, D_2,\dots, D_k, D_{k+1}$. Since the building is the limit of the holomorphic spheres intersecting the $J_\infty$-holomorphic hypersurface $\mathbb{CP}^{k-1}\times T^{2(n-k)}$ with intersection number $+1$. This means that exactly one of the components $D_1, D_2,\dots, D_m$, say $D_m$, of the building in the top level intersects this hypersurface; otherwise, the total intersection number will be greater than $+1$. Here, we use the fact that distinct holomorphic objects intersect positively. Therefore, the components $D_1, D_2,\dots, D_{m-1}$ are in the complement $ B^{2k}(r)\times T^{2(n-k)}$ of the hypersurface $\mathbb{CP}^{k-1}\times T^{2(n-k)}$. Compactifying $D_1, D_2,\dots, D_{m-1}$ to surfaces with boundaries on $L_{\Phi,t}$, we get
	 \[\frac{r}{k+1}\leq \frac{t}{k+1}\leq  \int D_j^*\Omega\]
	 for each $j=1,\dots, m-1$. Combining this with Equation (\ref{embed-est}), we have
	  \[\frac{r}{k+1}(m-1)+\int D_m^*\Omega\leq \sum_{j=1}^{m} \int D_j^*\Omega= r\]
	We must have $m-1\leq k+1$. Moreover, we cannot have $m-1=k+1$  because otherwise $\int D_m^*\Omega=0$ and this is not possible. After all, $D_m$ is non-constant and $J_\infty$-holomorphic.
	
	The conclusion is $m=k+1$, and since there are at least $k+1$ planes in the top level, the top level consists of the disks $D_1, D_2,\dots, D_k, D_{k+1}$ in which exactly one disk, say $D_{k+1}=: D_{\infty}$ intersects the hypersurface $\mathbb{CP}^{k-1}\times T^{2(n-k)}$ with intersection number $+1$. Moreover, for all $j=1,\dots,k$ we  have
	\[\frac{t}{k+1}= \int D_j^*\Omega\]
	 and 
	 \begin{equation}\label{engerydisk}
	 	0<\int D_\infty^*\Omega=r-\frac{t}{k+1}k.
	 \end{equation}

	Here, we want to point out that if we have the symplectic embedding in (\ref{symplecticembedding}), then from Equation \ref{engerydisk} we have $0<r-\frac{t}{k+1}k$ for every $t\in [r,k+1]$. In particular, for $t=k+1$ we have the obstruction $r>k$. However, the obstruction $r\geq k$ can also be derived from the  $k$-th Ekland--Hofer capacity.
\item[\textbf{Step 06}] Following the arguments of {\cite[Section 6.3]{Faisal:2024}}, we conclude that the holomorphic building has only two levels. The top level $\mathbb{CP}^k\times T^{2(n-k)}\setminus L_{\Phi,t}$ consists of $k+1$ negatively asymptotical cylindrical simple planes $D_1,\dots,D_k,D_\infty$. The bottom level $T^* L_{\Phi,t}$ consists of a single smooth connected punctured sphere $C_\mathrm{bot}$ with $k+1$ positive punctures carrying the tangency constraint at $p$. Moreover, each of the disks $D_1,\dots,D_k,D_\infty$ is of Maslov index $2$.
  
\item[\textbf{Step 07}] Let $\gamma_1,\gamma_2,\dots,\gamma_k, \gamma_{\infty}$ be the asymptotic closed Reeb orbits of the $J_\infty$-holomorphic planes $D_1,D_2,\dots,D_k, D_{\infty}$, respectively. Define	

\[
	\mathcal{M}^{J_\infty}_{[D_i]}(\gamma_i):=\left\{
	\begin{array}{l}
		D:\mathbb{CP}^1\setminus \{\infty\} \to (\mathbb{CP}^k\times T^{2(n-k)}\setminus L_{\Phi,t},J_{\infty}),\\
		du\circ i=J_{\infty}\circ du  ,\\
		u \text{ is asymptotic to  $\gamma_i$ at $\infty$,} \\
		
		\text{and } [D]=[D_i]\in \pi_2(\mathbb{CP}^k\times T^{2(n-k)},L_{\Phi,t}).
	\end{array}
	\right\}\bigg/\operatorname{Aut}(\mathbb{CP}^1,\infty)
\]
where $i=1,2,\dots,k,\infty$. By Equation (\ref{rel}) and Lemma \ref{count-imp}, we have 
	\[\# \mathcal{M}^{J_\infty}_{[D_i]}(\gamma_i)\neq 0\]
for all $i= 1, \dots, k, \infty$.
	
\item [\textbf{Step 08}]  For each $i= 1, \dots, k, \infty$, following {\cite[Section 6.6]{Faisal:2024}}, one can glue a half-cylinder to the moduli space $\mathcal{M}^{J_\infty}_{[D_i]}(\gamma_i)$ to create a Maslov index $2$  disk $\bar{D}_i:(D^2,\partial D^2)\to (\mathbb{CP}^{k}\times T^{2(n-k)}, L_{\Phi,t})$ whose boundary passes through a fixed generic point $q\in L_{\Phi,t}$. Moreover, this disk is $J$-holomorphic for generic $\Omega$-compatible almost complex structure $J$ on $\mathbb{CP}^{k}\times T^{2(n-k)}$. Let $\mathcal{M}^{J}_{L_{\Phi,t},[\bar{D}_i]}(q)$ denote the connected component of the moduli space containing  $\bar{D}_i$. By construction, we have the non-vanishing of the signed count 
\[\# \mathcal{M}^{J}_{L_{\Phi,t},[\bar{D}_i]}(q)\ne 0\]
for all $i= 1, \dots, k, \infty$.
Note that, by Equation (\ref{engerydisk}), every $D\in  \mathcal{M}^{J}_{L_{\Phi,t},[\bar{D}_\infty]}(q)$ has the symplectic area
\[0<\int D^*\Omega=r-\frac{t}{k+1}k.\]
This symplectic area is minimal in the sense that any non-trivial breaking in $\mathcal{M}^{J}_{L_{\Phi,t},[\bar{D}_\infty]}(q)$ requires a symplectic area strictly greater than this (cf.  {\cite[Theorem 6.10]{Faisal:2024}}). This implies the signed count above does not change if we vary $J$, $q$, or $t$ in the interval $[r,k+1]$. The same is true for other the moduli spaces $\mathcal{M}^{J}_{L_{\Phi,t},[\bar{D}_i]}(q)$. The conclusion is that the disks $\bar{D}_1,\bar{D}_2,\dots,\bar{D}_k, \bar{D}_{\infty}$ contribute non-trivially to the Landau--Ginzburg superpotential of $L_{\Phi,r}$ viewed as a Lagrangian torus in $(\mathbb{CP}^k\times\mathbb{C}^{n-k},r\omega_{\mathrm{FS}}\oplus\omega_{\mathrm{std}})$. Whereas, the disks $\bar{D}_1,\bar{D}_2,\dots,\bar{D}_k$  contribute to the Landau--Ginzburg superpotential of $L_{\Phi,r}$ viewed as a Lagrangian torus in $(B^{2k}(r)\times\mathbb{C}^{n-k},\omega_{\mathrm{FS}}\oplus\omega_{\mathrm{std}})$.

\item [\textbf{Step 09}] Our aim is to prove that $L_{\Phi,r}$ is not Hamiltonian isotopic to the Clifford torus. Define $e_i:=(0,\dots,1,\dots, 0)\in H_1(L_{\Phi,r}, \mathbb{Z})=\mathbb{Z}^n$ to be the class corresponding the $i$-th $S^1$-factor in $L_{\Phi,r}$. Suppose, on the contrary, that $L_{\Phi,r}$ is Hamiltonian isotopic to the Clifford torus, then $L_{\Phi,r}$ has Landau--Ginzburg potential of the Clifford torus. This means the disks $\bar{D}_1,\bar{D}_2,\dots,\bar{D}_k$ contribute to the Landau--Ginzburg potential of the Clifford torus in $(\mathbb{CP}^k\times\mathbb{C}^{n-k},r\omega_{\mathrm{FS}}\oplus\omega_{\mathrm{std}})$. So we must have that  $[\partial \bar{D}_i]=e_{l_i}\in H_1(L_{\Phi,r}, \mathbb{Z})=\mathbb{Z}^n$ for distinct $l_i\in \{1,\dots, n\}$, where $i=1,\dots,k$. The curve $C_\mathrm{bot}$ from \textbf{Step 06} provides a null-homology of the class $[\partial \bar{D}_\infty]+\sum_{i=1}^{k}[\partial \bar{D}_i]$, therefore, $[\partial \bar{D}_\infty]+\sum_{i=1}^{k}[\partial \bar{D}_i]=0$. This means $[\partial \bar{D}_\infty]=-\sum_{i=1}^{k}e_{l_i}$. 

\item [\textbf{Step 10}] 
Fix $t\in [r,k+1]$, choose $S>0$ such that 
	\[\frac{t}{k+1}<\frac{1}{1+(n-1)/S}<1.\]
The embedding $\Phi$ restricts to a symplectic embedding (still denoted by  $\Phi$)
	\[\Phi:E^{2n}(1,S,\dots, S)\to B^{2k}(r)\times \mathbb{C}^{n-k}.\]
This gives a symplectic embedding
	\[\Phi: E^{2n}(1,S,\dots, S)\to (\mathbb{CP}^{k}\times \mathbb{C}^{n-1}, r\omega_{\mathrm{FS}}\oplus \omega_{\mathrm{std}})\]
and the image of this embedding lies in the complement of the hypersurface $\mathbb{CP}^{k-1}\times \mathbb{C}^{n-k}$ for every large $S$. We have
\[L_{\Phi,t}:=\Phi\bigg(S^1\bigg(\frac{t}{k+1}\bigg)\times\cdots\times S^1\bigg(\frac{t}{k+1}\bigg)\bigg)\subset \Phi\big(\operatorname{int}(E^{2n}(1,S,\dots, S))\big)\]
for every large $S$. This allow us to apply neck-stretching to the moduli space of Maslov $2$ disks $\mathcal{M}^{J}_{L_{\Phi,t},[\bar{D}_\infty]}(q)$ along the contact type hypersurface $\partial E^{2n}(1,S,\dots, S)$ in $\mathbb{CP}^k\times T^{2(n-k)}$.

\item[\textbf{Step 11}] Take a sequence of almost complex structures $J_m$ on $\mathbb{CP}^k\times T^{2(n-k)}$ that stretches along the contact type hypersurface $\partial E^{2n}(1,S,\dots, S)$ in $\mathbb{CP}^k\times T^{2(n-k)}$. We assume $J_m$ restricted to a small neighborhood of the hypersurface $\mathbb{CP}^{k-1}\times T^{2(n-k)}$ is the standard complex structure $J_{\mathrm{std}}\oplus J_{\mathrm{std}}$ so that $\mathbb{CP}^{k-1}\times T^{2(n-k)}$ is $J_\infty$-holomorphic. Here $J_\infty$ is the almost complex structure on the symplectic completion of $\mathbb{CP}^k\times T^{2(n-k)}\setminus E^{2n}(1,S,\dots, S)$ obtained as the limit of $J_m$. Let $J_{\mathrm{bot}}$ denote the limiting SFT-admissible almost complex structure on the symplectic completion $\widehat{E}^{2n}(1,S,\dots, S)$.

\item[\textbf{Step 12}] Choose a sequence $D^m \in \mathcal{M}^{J_m}_{L_{\Phi,t},[\bar{D}_\infty]}(q)$. As $m\to \infty$, the $J_m$-holomorphic disk $D^m$ breaks to a holomorphic building $\mathbb{H}_\infty$ with the top level in the symplectic completion of $\mathbb{CP}^k\times T^{2(n-k)}\setminus E^{2n}(1,S,\dots, S)$, denoted by $\widehat{W}$, the bottom level in $\widehat{E^{2n}}(1,S,\dots, S)$ and some symplectization levels $\mathbb{R}\times \partial E^{2n}(1,S,\dots, S)$. We show that $\mathbb{H}_\infty$  consists of only two levels. The top level consists of a single degree one rigid $J_\infty$-holomorphic plane, denoted by $u_\infty$, asymptotic to the $k$-fold cover of the short Reeb orbit, denoted by $\beta^k$. The bottom level $\widehat{E^{2n}}(1,S,\dots, S)$ consists of a half-cylinder, denoted by $u_{\mathrm{cyl}}$, with boundary on $L_{\Phi,t}$ and positively asymptotic to $\beta^k$. See Figure \ref{embed01} for an illustration.
\begin{figure}[h]
\centering
\includegraphics[width=8cm]{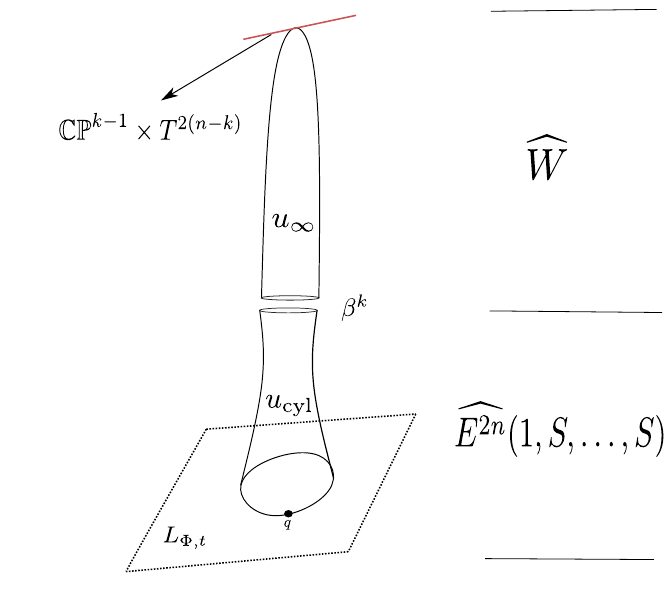}
\caption{$\widehat{W}$ denotes the symplectic completion of $\mathbb{CP}^k\times T^{2(n-k)}\setminus E^{2n}(1,S,\dots, S)$.}\label{embed01}
\end{figure}

\item[\textbf{Step 13}]  Let $u_1, u_2,\dots, u_l, u_\infty$ denote the smooth connected components of the building $\mathbb{H}_\infty$ in the top level. By construction, exactly one curve component, say $u_\infty$, intersects the complex hypersurface $\mathbb{CP}^{(k-1)}\times T^{2(n-k)}$, and the intersection number is $+1$. The other components $u_1, u_2,\dots, u_l$ are contained in the complement of $\mathbb{CP}^{(k-1)}\times T^{2(n-k)}$. Moreover, none of these is a closed $J_\infty$-holomorphic sphere because the symplectic form is exact on the complement of $\mathbb{CP}^{(k-1)}\times T^{2(n-k)}$. Because the intersection number of $u_\infty$ with $\mathbb{CP}^{(k-1)}\times T^{2(n-k)}$ is $+1$, $u_\infty$ somewhere injective. Moreover, the bottom level $\widehat{E^{2n}}(1,S,\dots, S)$ contains a smooth connected curve, denoted by $u_{\mathrm{cyl}}$,  with some positive punctures asymptotic to closed Reeb orbits on $\partial E^{2n}(1,S,\dots, S)$ and boundary on $L_{\Phi,t}$. The boundary of $u_{\mathrm{cyl}}$ represents the class $[\partial \bar{D}_\infty]=-\sum_{i=1}^{k}e_{l_i}\in H_1(L_{\Phi,r}, \mathbb{Z})$ by \textbf{Step 09}.

\item[\textbf{Step 14}] The energy of $\mathbb{H}_\infty$ is given by Equation (\ref{engerydisk}). We have
\begin{equation}\label{symengerestem}
0<\int u_\infty^*\tilde{\Omega}+\sum_{i=1}^{l}\int u_i^*\tilde{\Omega}\leq  r-\frac{t}{k+1}k<1.
\end{equation}
Here,
\begin{equation*}
	\tilde{\Omega}:=
	\begin{cases}
		r\omega_{\mathrm{FS}}\oplus \omega_{\mathrm{std}}& \text{on } \mathbb{CP}^k\times T^{2(n-k)}\setminus E^{2n}(1,S,\dots, S),\\
		d\lambda_{\mathrm{std}}|_{\partial E^{2n}(1,S,\dots, S)} & \text{on } (-\infty,0]\times \partial E^{2n}(1,S,\dots, S).
	\end{cases}
\end{equation*}
The curve components $u_1, u_2,\dots, u_l$ are contained in the complement of $\mathbb{CP}^{(k-1)}\times T^{2(n-k)}$ where the symplectic form is exact, and are negatively asymptotic to closed Reeb orbits on $\partial E^{2n}(1,S,\dots, S)$. By Section \ref{dynomicsofelipsoids}, the minimal period of a closed Reeb orbit on $\partial E^{2n}(1,S,\dots, S)$ is $1$. Thus,
\[\sum_{i=1}^{l}\int u_i^*\tilde{\Omega}\geq l.\]
This contradicts (\ref{symengerestem}) if $l\geq 1$. This means $l=0$, i.e.,  curve components $u_1, u_2,\dots, u_l$ do not exist. So the top level consists of a single simple smooth connected curve $u_\infty$.
\item[\textbf{Step 15}]
We prove that all negative punctures of $u_\infty$ are asymptotic to covers (possibly multiple) of the short Reeb orbit $\beta_1$. Suppose  $u_\infty$ has negative ends on the Reeb orbits $\beta^{m_1}_{i_1}, \beta^{m_2}_{i_2},\dots,\beta^{m_l}_{i_l}$ and assume at least one, say $\beta^{m_1}_{i_1}$, is a long orbit. The Fredholm index of $u_\infty$ in the trivialization $\tau_{\mathrm{ext}}$ is
\[\operatorname{ind}(u_\infty)=(n-3)(2-l)+2c _1([\mathbb{CP}^1])-\sum_{j=1}^{l}\operatorname{CZ}^{\tau_{\mathrm{ext}}}(\beta^{m_j}_{i_j}).\] 
One can see that
\begin{equation}\label{local1}
\operatorname{ind}(u_\infty)\leq (n-3)+2c _1([\mathbb{CP}^1])-\operatorname{CZ}^{\tau_{\mathrm{ext}}}(\beta^{m_1}_{i_1})..
\end{equation}
By Theorem \ref{CZIimp}, for the long orbit $\beta^{m_1}_{i_1}$ we have
\begin{equation}\label{local2}
\operatorname{CZ}^{\tau_{\mathrm{ext}}}(\beta^{m_1}_{i_1})\geq n-1+2(\lfloor S\rfloor+n-1).
\end{equation}

Combining Equations (\ref{local1}) and (\ref{local2}) yields
	\[\operatorname{ind}(u_\infty)\leq 2(c _1([\mathbb{CP}^1])-\lfloor S\rfloor-n).\] 
	Note that  $c _1([\mathbb{CP}^1])=k+1$ and also by our assumption $S\geq k$, therefore
	\[\operatorname{ind}(u_\infty)\leq 2(k-1-\lfloor S\rfloor) \leq -2. \]
The curve $u_\infty$ is simple and  $J_\infty$-holomorphic. We can perturb $J_\infty$ near the hypersurface $\mathbb{CP}^{k-1}\times T^{2(n-k)}$ to assume $u_\infty$ is regular. So we must have $\operatorname{ind}(u_\infty)\geq 0$. This is a contradiction to the above estimate. Thus, all the ends of  $u_\infty$ are on short Reeb orbits.

\item[\textbf{Step 16}]  Next we prove that $u_\infty$ has a single negative puncture that is asymptotic to the $m$-fold cover of the short Reeb orbit, denoted by $\beta^m$, for some positive integer $m\leq k$. By \textbf{Step 15}, all negative punctures of $u_\infty$ are asymptotic to short  Reeb orbits. Suppose  $u_\infty$ has negative ends on the Reeb orbits $\beta^{m_1}, \beta^{m_2},\dots,\beta^{m_l}$. The Fredholm index of $u_\infty$ in the trivialization $\tau_{\mathrm{ext}}$ is
\[\operatorname{ind}(u_\infty)=(n-3)(2-l)+2(k+1)-\sum_{j=1}^{l}\operatorname{CZ}^{\tau_{\mathrm{ext}}}(\beta^{m_j}).\] 
By Theorem \ref{CZIimp}, we have
\[\operatorname{CZ}^{\tau_{\mathrm{ext}}}(\beta^{m_i})\geq n-1+2m_i. \]
This implies
\[\operatorname{ind}(u_\infty)\leq (n-3)(2-l)+2(k+1)-l(n-1)-2\sum_{i=1}^{l}m_i.\]
If $l\geq 2$, then 
\[	\operatorname{ind}(u_\infty)\leq 2(k+1)-2(n-1)-2\sum_{i=1}^{2}m_i\]
Since $n\geq k+1$ and $\sum_{i=1}^{2}m_i\geq 2$, we have 
\[\operatorname{ind}(u_\infty)\leq -2.\]
This is again a contradiction. So we must have $l=1$, i.e., $u_\infty$ has only one negative puncture.

Suppose the negative puncture of $u_\infty$ is asymptotic to $\beta^m$. By the same arguments as above, we have 
\[0\leq \operatorname{ind}(u_\infty)\leq 2(k-m).\]
This means $m\leq k$.

\item[\textbf{Step 17}] Recall from \textbf{Step 13} that the bottom level $\widehat{E^{2n}}(1,S,\dots, S)$ contains a smooth connected curve with connected boundary, denoted by $u_{\mathrm{cyl}}$,  with some positive punctures asymptotic to closed Reeb orbits on $\partial E^{2n}(1,S,\dots, S)$ and boundary on $L_{\Phi,t}$. The boundary of $u_{\mathrm{cyl}}$ represents the class $[\partial \bar{D}_\infty]=-\sum_{i=1}^{k}e_{l_i}\in H_1(L_{\Phi,r}, \mathbb{Z})$. We prove that $u_{\mathrm{cyl}}$ has only one postive puncture and is therefore a half-cylinder $u_{\mathrm{cyl}}:[0,\infty)\times S^1 \to \widehat{E^{2n}}(1,S,\dots, S) $ with $u({\{0\}\times S^1})\subset L_{\Phi,t}$.

The underlying graph of the building $\mathbb{H}_\infty$ is a tree since the building has genus zero. Suppose  $u_{\mathrm{cyl}}$ has $m$ positive punctures, for some positive integer $m$. There are $m$ edges emanating from the vertex $u_{\mathrm{cyl}}$ in the underlying graph. We order these edges from $1,2, \dots, m$. Let $C_i$ be the subtree emanating from the vertex $u_{\mathrm{cyl}}$ along the $i$th edge. The trees $C_1,\dots, C_{k+1},\dots, C_{m}$ are topological planes with curve components in different levels.  Since the building has only one curve component in the top level, and that is  $u_\infty$, at most one of $C_1,\dots, C_{m}$, say $C_{m}$, contains $u_\infty$. By the maximum principle, each of $C_i$ must have some curve components in the top level.  Thus, we have at least $m$ smooth connected components in the top level. But by \textbf{Step 14}, there is only one curve component in the top level, namely $u_\infty$. Thus, we must have $m=1$,  i.e.,  $u_{\mathrm{cyl}}$ has only one positive puncture. 

\item[\textbf{Step 18}] Next, we prove that the positive puncture of $u_{\mathrm{cyl}}$ is asymptotic to $\beta^k$, the $k$-fold cover of the short Reeb orbit $\beta_1$. Suppose the positive puncture of $u_{\mathrm{cyl}}$ is asymptotic to $\beta^l$, for some positive integer $l$. The boundary of $u_{\mathrm{cyl}}$ on $L_{\Phi,t}$ bounds a disk of symplectic area $tk/(k+1)$, so 
\[0\leq \int u_{\mathrm{cyl}}^*\tilde{\omega}_{\mathrm{std}}=\int_{\beta^l}\lambda_{\mathrm{std}}-\frac{t\, k}{k+1} \]
where $\tilde{\omega}_{\mathrm{std}}$ is the exact $2$-form
\begin{equation*}
	\tilde{\omega}_{\mathrm{std}}:=
	\begin{cases}
		d\lambda_{\mathrm{std}}|_{\partial E^{2n}(1,S,\dots, S)}& \text{on } [0,\infty)\times \partial E^{2n}(1,S,\dots, S),\\
		\omega_{\mathrm{std}} & \text{on } E^{2n}(1,S,\dots, S)\setminus L_{\Phi,t}.
	\end{cases}
\end{equation*}
Thus, we have 
\[\frac{t\, k}{k+1}\leq\int_{\beta^l}\lambda_{\mathrm{std}}=l. \]
Since $l$ is an integer and $t$ can be chosen arbitrary closed to $k+1$, therefore  $l\geq k.$

Recall that $u_\infty$ is asymptotic to $\beta^m$ for some $m\leq k$. The total action of Reeb orbits decreases as one goes from top to bottom along the building in the symplectization levels  $\mathbb{R}\times  \partial E^{2n}(1,S,\dots, S)$, so we must have $k\geq m\geq l\geq k$. Putting everything together, we obtain $m=l=k$.

The conclusion is that $u_\infty$ is negatively asymptotic to $\beta^k$ and $u_{\mathrm{cyl}}$ is positively asymptotic to $\beta^k$. 

\item[\textbf{Step 19}] Recall that $u_\infty$ is somewhere injective. Also, $u_{\mathrm{cyl}}$ is somewhere injective because its boundary represents a primitive class in $H_1(L_{\Phi,t}, \mathbb{Z})$. So, generically, both $u_{\mathrm{cyl}}$ and $u_\infty$ have non-negative Fredholm indices. The index of the building $\mathbb{H}_\infty$ is zero, so 
\[\underbrace{\operatorname{ind}(u_{\mathrm{cyl}})}_{\geq 0}+\underbrace{\operatorname{ind}(u_\infty)}_{\geq 0}=0\]
which implies $\operatorname{ind}(u_{\mathrm{cyl}})=\operatorname{ind}(u_\infty)=0$.

\item[\textbf{Step 20}] The $\tilde{\omega}_{\mathrm{std}}$-energy of the rigid half-cylinder $u_{\mathrm{cyl}}:[0,\infty)\times S^1 \to \widehat{E^{2n}}(1,S,\dots, S) $  satisfies 
\[0\leq \int u_{\mathrm{cyl}}^*\tilde{\omega}_{\mathrm{std}}= \int_{u_{\mathrm{cyl}}^{-1}(E^{2n}(1,S,\dots, S))}u_{\mathrm{cyl}}^*\omega_{\mathrm{std}}+\int_{u_{\mathrm{cyl}}^{-1}([0,\infty)\times \partial E^{2n}(1,S,\dots, S))} u_{\mathrm{cyl}}^*d\lambda_{\mathrm{std}}=k-\frac{t\, k}{k+1},\]
where  $t\in [r,k+1)$. So the $\tilde{\omega}_{\mathrm{std}}$-energy of $u_{\mathrm{cyl}}$ can be made arbitrary small by moving $t$ toward $k+1$. In case $t=k+1$, we have 
\[0=\int_{u_{\mathrm{cyl}}^{-1}(E^{2n}(1,\infty,\dots, \infty))}u_{\mathrm{cyl}}^*\omega_{\mathrm{std}}+\int_{u_{\mathrm{cyl}}^{-1}([0,\infty)\times \partial E^{2n}(1,\infty,\dots, \infty))} u_{\mathrm{cyl}}^*d\lambda_{\mathrm{std}}.\]
This implies 
\begin{equation}\label{last1}
\int_{u_{\mathrm{cyl}}^{-1}(E^{2n}(1,\infty,\dots, \infty))}u_{\mathrm{cyl}}^*\omega_{\mathrm{std}}=0
\end{equation}
and 
\begin{equation}\label{last2}
	\int_{u_{\mathrm{cyl}}^{-1}([0,\infty)\times \partial E^{2n}(1,\infty,\dots, \infty))} u_{\mathrm{cyl}}^*d\lambda_{\mathrm{std}}=0.
\end{equation}
Equation \ref{last1} implies $u_{\mathrm{cyl}}$ is entirely contained in the cylindrical
 end, i.e., 
\[u_{\mathrm{cyl}}:[0,\infty)\times S^1 \to [0,\infty)\times \partial E^{2n}(1,\infty,\dots, \infty)\]
 with $u_{\mathrm{cyl}}(\{0\}\times S^1)\subset L_{\Phi,k+1}\subset  \{0\}\times \partial E^{2n}(1,\infty,\dots, \infty)$. Note that $L_{\Phi,k+1}=S^1(1)\times \cdots\times S^1(1).$

Next, we explain an implication of Equation \ref{last2}. Let $\lambda_{\mathrm{std}}:=\lambda_{\mathrm{std}}|_{\partial E^{2n}(1,\infty,\dots, \infty)}$ and set $\xi_{\mathrm{std}}:=\operatorname{Ker}(\lambda_{\mathrm{std}}).$ Denote by $R_{\mathrm{std}}$  the Reeb vector field of $\lambda_{\mathrm{std}}$. We have the splitting 
\[T([0,\infty)\times \partial E^{2n}(1,\infty,\dots,\infty))=\operatorname{span}\{\partial_{r}, R_{\mathrm{std}}\}\oplus \xi_{\mathrm{std}},\]
where $\partial_{r}$ is the unit vector field in the $\mathbb{R}$-direction. The $2$-form $d\lambda_{\mathrm{std}}|_{\xi_{\mathrm{std}}}$ is symplectic  and $J_{\mathrm{bot}}$ is $d\lambda_{\mathrm{std}}$-compatible on $\xi_{\mathrm{std}}$. Moreover, $d\lambda_{\mathrm{std}}$ annihilates the trivial subbundle $\operatorname{span}\{\partial_{r}, R_{\mathrm{std}}\}$. Let $\pi_{\xi_{\mathrm{std}}}:T([0,\infty)\times \partial E^{2n}(1,\infty,\dots,\infty))\to \xi_{\mathrm{std}}$
be the projection along the trivial subbundle $\operatorname{span}\{\partial_{r}, R_{\mathrm{std}}\}$. Since $u_{\mathrm{cyl}}$ is $J_{\mathrm{bot}}$-holomorphic, for any $(s,t)\in [0,\infty)\times S^1$, we have
\begin{equation}\label{sourceofcontra}
	d\lambda_{\mathrm{std}}(\partial _s u_{\mathrm{cyl}}, \partial _t u_{\mathrm{cyl}})=d\lambda_{\mathrm{std}}(\partial _s u_{\mathrm{cyl}}, J_{\mathrm{bot}}\partial _s u_{\mathrm{cyl}})=d\lambda_{\mathrm{std}}(\pi_{\xi_{\mathrm{std}}}\partial _t u_{\mathrm{cyl}}, J_{\mathrm{bot}}\pi_{\xi_{\mathrm{std}}}\partial _t u_{\mathrm{cyl}})\geq 0,
\end{equation}
where the inequality is strict if $\pi_{\xi_{\mathrm{std}}}\partial _t u_{\mathrm{cyl}}$ does not vanish.

The loop $u_{\mathrm{cyl}}(\{0\}\times S^1)\subset L_{\Phi,k+1}\subset \{0\}\times \partial E^{2n}(1,\infty,\dots, \infty)$, where $L_{\Phi,k+1}=S^1(1)\times \cdots\times S^1(1)$, parameterizes a closed Reeb orbit. Suppose, on the contrary, that it is not the case. Then for some $t_0\in S^1$, we must have $\pi_{\xi_{\mathrm{std}}}\partial _t u_{\mathrm{cyl}}(0,t_0)\neq 0$. By continuity, for some $\epsilon_1,\epsilon_2>0$, we have $\pi_{\xi_{\mathrm{std}}}\partial _t u_{\mathrm{cyl}}(s,t_0)\neq 0$ for all $(s,t)\in [0,\epsilon_1)\times (t_0-\epsilon_2,t_0+\epsilon_2)$. By (\ref{sourceofcontra}) and Equation (\ref{last2}), we obtain the contradiction
\[0= \int_{u_{\mathrm{cyl}}^{-1}([0,\infty)\times \partial E^{2n}(1,\infty,\dots, \infty))} u_{\mathrm{cyl}}^*d\lambda_{\mathrm{std}}\geq \int_{u_{\mathrm{cyl}}^{-1}([0,\epsilon_1)\times (t_0-\epsilon_2,t_0+\epsilon_2))} u_{\mathrm{cyl}}^*d\lambda_{\mathrm{std}}>0.\]

The standard contact form $\lambda_{\mathrm{std}}$ on  $\partial E^{2n}(1,\infty,\dots, \infty)$ is Morse--Bott. Simple closed Reeb orbits are of the form $\beta^1\times\{z\}$ for $z\in \mathbb{C}^{n-1}$ and vice versa. In particular, there is no closed Reeb orbit on $L_{\Phi,k+1}=S^1(1)\times \cdots\times S^1(1)$ in the primitive class \[[u_{\mathrm{cyl}}(\{0\}\times S^1)]=[\partial \bar{D}_\infty]=-\sum_{i=1}^{k}e_{l_i}\in H_1(L_{\Phi,k+1}, \mathbb{Z}).\] But by the argument above the loop  $u_{\mathrm{cyl}}(\{0\}\times S^1)\subset L_{\Phi,k+1}\subset \partial E^{2n}(1,\infty,\dots, \infty)$, where $L_{\Phi,k+1}=S^1(1)\times \cdots\times S^1(1)$, parameterizes a closed Reeb orbit--- which means we must have that $[\partial \bar{D}_\infty]=-ke_{1}\in H_1(L_{\Phi,k+1}, \mathbb{Z})$. This is a contradiction. This completes our proof.
\end{itemize}


\AtEndDocument{
\bibliographystyle{alpha}
\bibliography{extremallagrangian}

@article{Chaidez:2024aa,
	author = {Julian Chaidez},
	eprint = {2412.20676},
	journal = {International Mathematics Research Notices},
	month = {12},
	number = {21},
	title = {Contact homology and linearization without dga homotopies},
	url = {https://arxiv.org/pdf/2412.20676.pdf},
	volume = {2025},
	year = {2025}}

@article{Buse-Hindellipsoidsindimensiongreaterthanfour,
	author = {Buse, Olguta and Hind, Richard},
	doi = {10.2140/gt.2011.15.2091},
	fjournal = {Geometry \& Topology},
	issn = {1465-3060},
	journal = {Geom. Topol.},
	mrclass = {53D35},
	mrnumber = {2860988},
	mrreviewer = {Emmanuel Opshtein},
	number = {4},
	pages = {2091--2110},
	title = {Symplectic embeddings of ellipsoids in dimension greater than four},
	url = {https://doi.org/10.2140/gt.2011.15.2091},
	volume = {15},
	year = {2011}}

@article{McDuff:2024ab,
	author = {Dusa McDuff and Kyler Siegel},
	eprint = {2412.00561},
	journal = {arXiv:2412.00561v2},
	month = {12},
	title = {Sesquicuspidal curves, scattering diagrams, and symplectic nonsqueezing},
	url = {https://arxiv.org/pdf/2412.00561.pdf},
	year = {2024}}

@article{Pereira:2025aa,
	author = {Pereira, Miguel},
	doi = {10.4310/jsg.250403022006},
	fjournal = {The Journal of Symplectic Geometry},
	issn = {1527-5256},
	journal = {J. Symplectic Geom.},
	mrclass = {53D05},
	mrnumber = {4886504},
	number = {1},
	pages = {159--225},
	title = {On the {L}agrangian capacity of convex or concave toric domains},
	url = {https://doi.org/10.4310/jsg.250403022006},
	volume = {23},
	year = {2025}}

@article{Schlenkellipsoids,
	author = {Schlenk, Felix},
	doi = {10.1007/BF02783427},
	fjournal = {Israel Journal of Mathematics},
	issn = {0021-2172},
	journal = {Israel J. Math.},
	mrclass = {53D35 (57R17)},
	mrnumber = {2031958},
	mrreviewer = {Amine Hadjar},
	pages = {215--252},
	title = {Symplectic embeddings of ellipsoids},
	url = {https://doi.org/10.1007/BF02783427},
	volume = {138},
	year = {2003}}

@article{Faisal:2024,
	author = {Shah Faisal},
	journal = {arXiv:2504.13076},
	pages = {91},
	title = {Extremal {L}agrangian tori in toric domains},
	year = {2025}}

@article{Faisal:2024aa,
	author = {Shah Faisal},
	eprint = {2412.18462},
	journal = {arXiv:2412.18462, 2024},
	title = {A proof of {G}romov's non-squeezing theorem},
	url = {https://arxiv.org/pdf/2412.18462.pdf}}

@article{Wendl:2016aa,
	author = {Chris Wendl},
	eprint = {1612.01009},
	journal = {arXiv:1612.01009},
	month = {12},
	title = {Lectures on Symplectic Field Theory},
	url = {https://arxiv.org/pdf/1612.01009.pdf},
	year = {2016}}

@article{Pardon-Contacthomologyandvirtualfundamentalcycles,
	author = {Pardon, John},
	doi = {10.1090/jams/924},
	fjournal = {Journal of the American Mathematical Society},
	issn = {0894-0347},
	journal = {J. Amer. Math. Soc.},
	mrclass = {53D42 (53D35 53D40 57R17)},
	mrnumber = {3981989},
	mrreviewer = {Alexander Fel\cprime shtyn},
	number = {3},
	pages = {825--919},
	title = {Contact homology and virtual fundamental cycles},
	url = {https://doi.org/10.1090/jams/924},
	volume = {32},
	year = {2019}}

@article{McDuff:2023aa,
	author = {Dusa McDuff and Kyler Siegel},
	eprint = {2308.07542},
	journal = {arXiv:2308.07542},
	month = {08},
	title = {Ellipsoidal superpotentials and singular curve counts},
	url = {https://arxiv.org/pdf/2308.07542.pdf},
	year = {2023}}

@article{McDuff:2024aa,
	author = {Dusa McDuff and Kyler Siegel},
	eprint = {2404.14702},
	journal = {arXiv:2404.14702},
	month = {04},
	title = {Singular algebraic curves and infinite symplectic staircases},
	url = {https://arxiv.org/pdf/2404.14702.pdf},
	year = {2024}}

@article{MR2026549,
	author = {Bourgeois, Fr\'{e}d\'{e}ric and Eliashberg, Yakov and Hofer, Helmut and Wysocki, Kris and Zehnder, Eduard},
	doi = {10.2140/gt.2003.7.799},
	fjournal = {Geometry and Topology},
	issn = {1465-3060},
	journal = {Geom. Topol.},
	mrclass = {53D45 (53D35 53D40 57R17)},
	mrnumber = {2026549},
	mrreviewer = {Kai Cieliebak},
	pages = {799--888},
	title = {Compactness results in symplectic field theory},
	url = {https://doi.org/10.2140/gt.2003.7.799},
	volume = {7},
	year = {2003}}

@article{McDuff-Hoferconjecture,
	author = {McDuff, Dusa},
	fjournal = {Journal of Differential Geometry},
	issn = {0022-040X},
	journal = {J. Differential Geom.},
	mrclass = {53D35},
	mrnumber = {2844441},
	mrreviewer = {Sheila Sandon},
	number = {3},
	pages = {519--532},
	title = {The {H}ofer conjecture on embedding symplectic ellipsoids},
	url = {http://projecteuclid.org/euclid.jdg/1321366358},
	volume = {88},
	year = {2011}}

@article{McDuff-Schlenk,
	author = {McDuff, Dusa and Schlenk, Felix},
	doi = {10.4007/annals.2012.175.3.5},
	fjournal = {Annals of Mathematics. Second Series},
	issn = {0003-486X},
	journal = {Ann. of Math. (2)},
	mrclass = {53D42 (11B39 53D35)},
	mrnumber = {2912705},
	mrreviewer = {Sheila Sandon},
	number = {3},
	pages = {1191--1282},
	title = {The embedding capacity of 4-dimensional symplectic ellipsoids},
	url = {https://doi.org/10.4007/annals.2012.175.3.5},
	volume = {175},
	year = {2012}}

@incollection{MR2369441,
	author = {Cieliebak, Kai and Hofer, Helmut and Latschev, Janko and Schlenk, Felix},
	booktitle = {Dynamics, ergodic theory, and geometry},
	doi = {10.1017/CBO9780511755187.002},
	mrclass = {53D35 (37J05 57R17)},
	mrnumber = {2369441},
	mrreviewer = {Michael J. Usher},
	pages = {1--44},
	publisher = {Cambridge Univ. Press, Cambridge},
	series = {Math. Sci. Res. Inst. Publ.},
	title = {Quantitative symplectic geometry},
	url = {https://doi.org/10.1017/CBO9780511755187.002},
	volume = {54},
	year = {2007}}

@article{MR3777016,
	author = {Schlenk, Felix},
	doi = {10.1090/bull/1587},
	fjournal = {American Mathematical Society. Bulletin. New Series},
	issn = {0273-0979},
	journal = {Bull. Amer. Math. Soc. (N.S.)},
	mrclass = {53D35 (37B40 53D40)},
	mrnumber = {3777016},
	mrreviewer = {Richard Keith Hind},
	number = {2},
	pages = {139--182},
	title = {Symplectic embedding problems, old and new},
	url = {https://doi.org/10.1090/bull/1587},
	volume = {55},
	year = {2018}}

@article{MR2475400,
	author = {Bourgeois, Fr\'{e}d\'{e}ric and Oancea, Alexandru},
	doi = {10.1215/00127094-2008-062},
	fjournal = {Duke Mathematical Journal},
	issn = {0012-7094},
	journal = {Duke Math. J.},
	mrclass = {53D40},
	mrnumber = {2475400},
	mrreviewer = {Tobias Ekholm},
	number = {1},
	pages = {71--174},
	title = {Symplectic homology, autonomous {H}amiltonians, and {M}orse-{B}ott moduli spaces},
	url = {https://doi.org/10.1215/00127094-2008-062},
	volume = {146},
	year = {2009}}

@article{MR3671507,
	author = {Bourgeois, Fr\'{e}d\'{e}ric and Oancea, Alexandru},
	doi = {10.1093/imrn/rnw029},
	fjournal = {International Mathematics Research Notices. IMRN},
	issn = {1073-7928},
	journal = {Int. Math. Res. Not. IMRN},
	mrclass = {53D42 (53D40 55N91 57R17 57R58)},
	mrnumber = {3671507},
	mrreviewer = {Daniel V. Mathews},
	number = {13},
	pages = {3849--3937},
	title = {{$S^1$}-equivariant symplectic homology and linearized contact homology},
	url = {https://doi.org/10.1093/imrn/rnw029},
	year = {2017}}

@article{MR3868228,
	author = {Gutt, Jean and Hutchings, Michael},
	doi = {10.2140/agt.2018.18.3537},
	fjournal = {Algebraic \& Geometric Topology},
	issn = {1472-2747},
	journal = {Algebr. Geom. Topol.},
	mrclass = {53D40 (53D05 57R17)},
	mrnumber = {3868228},
	mrreviewer = {Stephan Mescher},
	number = {6},
	pages = {3537--3600},
	title = {Symplectic capacities from positive {$S^1$}-equivariant symplectic homology},
	url = {https://doi.org/10.2140/agt.2018.18.3537},
	volume = {18},
	year = {2018}}

@article{Mikhalkin:2023aa,
	author = {Grigory Mikhalkin and Kyler Siegel},
	eprint = {2307.13252},
	journal = {arXiv:2307.13252},
	month = {07},
	title = {Ellipsoidal superpotentials and stationary descendants},
	url = {https://arxiv.org/pdf/2307.13252.pdf},
	year = {2023}}

@article{Cristofaro-Hind-Siegel-stabilizedembedding,
	author = {Cristofaro-Gardiner, Daniel and Hind, Richard and Siegel, Kyler},
	doi = {10.1007/s11784-022-00942-z},
	fjournal = {Journal of Fixed Point Theory and Applications},
	issn = {1661-7738},
	journal = {J. Fixed Point Theory Appl.},
	mrclass = {53D05},
	mrnumber = {4439980},
	mrreviewer = {Chun-Gen Liu},
	number = {2},
	pages = {Paper No. 49, 25},
	title = {Higher symplectic capacities and the stabilized embedding problem for integral elllipsoids},
	url = {https://doi.org/10.1007/s11784-022-00942-z},
	volume = {24},
	year = {2022}}

@article{MR4466005,
	author = {Siegel, Kyler},
	doi = {10.1093/imrn/rnaa334},
	fjournal = {International Mathematics Research Notices. IMRN},
	issn = {1073-7928},
	journal = {Int. Math. Res. Not. IMRN},
	mrclass = {53D05 (32Q65 53D12)},
	mrnumber = {4466005},
	mrreviewer = {Emmanuel Opshtein},
	number = {16},
	pages = {12402--12461},
	title = {Computing higher symplectic capacities {I}},
	url = {https://doi.org/10.1093/imrn/rnaa334},
	year = {2022}}

@article{MR3867635,
	author = {Hind, Richard and Kerman, Ely},
	doi = {10.1007/s00222-018-0828-7},
	fjournal = {Inventiones Mathematicae},
	issn = {0020-9910},
	journal = {Invent. Math.},
	mrclass = {53D05 (53C42)},
	mrnumber = {3867635},
	number = {2},
	pages = {1023--1029},
	title = {Correction to: {N}ew obstructions to symplectic embeddings},
	url = {https://doi.org/10.1007/s00222-018-0828-7},
	volume = {214},
	year = {2018}}

@article{MR3394319,
	author = {Hind, Richard},
	doi = {10.1112/jtopol/jtv016},
	fjournal = {Journal of Topology},
	issn = {1753-8416},
	journal = {J. Topol.},
	mrclass = {53D35 (57R17)},
	mrnumber = {3394319},
	mrreviewer = {Janko Latschev},
	number = {3},
	pages = {871--883},
	title = {Some optimal embeddings of symplectic ellipsoids},
	url = {https://doi.org/10.1112/jtopol/jtv016},
	volume = {8},
	year = {2015}}

@article{Hind-Kerman-ellipembedding,
	author = {Hind, Richard and Kerman, Ely},
	doi = {10.1007/s00222-013-0471-2},
	fjournal = {Inventiones Mathematicae},
	issn = {0020-9910},
	journal = {Invent. Math.},
	mrclass = {53D05 (53C42)},
	mrnumber = {3193752},
	mrreviewer = {Emmanuel Opshtein},
	number = {2},
	pages = {383--452},
	title = {New obstructions to symplectic embeddings},
	url = {https://doi.org/10.1007/s00222-013-0471-2},
	volume = {196},
	year = {2014}}

@article{MR3782228,
	author = {McDuff, Dusa},
	doi = {10.1007/s40879-017-0184-y},
	fjournal = {European Journal of Mathematics},
	issn = {2199-675X},
	journal = {Eur. J. Math.},
	mrclass = {53D05 (57R17)},
	mrnumber = {3782228},
	mrreviewer = {Michael J. Usher},
	number = {1},
	pages = {356--371},
	title = {A remark on the stabilized symplectic embedding problem for ellipsoids},
	url = {https://doi.org/10.1007/s40879-017-0184-y},
	volume = {4},
	year = {2018}}

@article{MR3789827,
	author = {Cristofaro-Gardiner, Daniel and Hind, Richard and McDuff, Dusa},
	doi = {10.1112/topo.12055},
	fjournal = {Journal of Topology},
	issn = {1753-8416},
	journal = {J. Topol.},
	mrclass = {57R17 (53D05)},
	mrnumber = {3789827},
	mrreviewer = {Emmanuel Opshtein},
	number = {2},
	pages = {309--378},
	title = {The ghost stairs stabilize to sharp symplectic embedding obstructions},
	url = {https://doi.org/10.1112/topo.12055},
	volume = {11},
	year = {2018}}

@article{MR4332489,
	author = {McDuff, Dusa and Siegel, Kyler},
	doi = {10.1112/topo.12204},
	fjournal = {Journal of Topology},
	issn = {1753-8416},
	journal = {J. Topol.},
	mrclass = {53D45},
	mrnumber = {4332489},
	mrreviewer = {Pierrick Bousseau},
	number = {4},
	pages = {1176--1242},
	title = {Counting curves with local tangency constraints},
	url = {https://doi.org/10.1112/topo.12204},
	volume = {14},
	year = {2021}}

@article{Hind,
	author = {Hind, Richard and Opshtein, Emmanuel},
	doi = {10.4171/CMH/496},
	fjournal = {Commentarii Mathematici Helvetici. A Journal of the Swiss Mathematical Society},
	issn = {0010-2571},
	journal = {Comment. Math. Helv.},
	mrclass = {53D12 (53D35 53E50)},
	mrnumber = {4152624},
	mrreviewer = {W.-H. Steeb},
	number = {3},
	pages = {535--567},
	title = {Squeezing {L}agrangian tori in dimension 4},
	url = {https://doi.org/10.4171/CMH/496},
	volume = {95},
	year = {2020}}

@article{Cieliebak_2005,
	author = {Kai Cieliebak and Klaus Mohnke},
	doi = {10.4310/jsg.2005.v3.n4.a5},
	journal = {Journal of Symplectic Geometry},
	number = {4},
	pages = {589--654},
	publisher = {International Press of Boston},
	title = {Compactness for punctured holomorphic curves},
	url = {https://doi.org/10.4310%2Fjsg.2005.v3.n4.a5},
	volume = {3},
	year = 2005}

@article{Cieliebak_2007,
	author = {Kai Cieliebak and Klaus Mohnke},
	doi = {10.4310/jsg.2007.v5.n3.a2},
	journal = {Journal of Symplectic Geometry},
	number = {3},
	pages = {281--356},
	publisher = {International Press of Boston},
	title = {Symplectic hypersurfaces and transversality in Gromov--Witten theory},
	url = {https://doi.org/10.4310%2Fjsg.2007.v5.n3.a2},
	volume = {5},
	year = 2007}

@book{Frauenfelder_2018,
	author = {Urs Frauenfelder and Otto van Koert},
	doi = {10.1007/978-3-319-72278-8},
	publisher = {Springer International Publishing},
	title = {The Restricted Three-Body Problem and Holomorphic Curves},
	url = {https://doi.org/10.1007%2F978-3-319-72278-8},
	year = 2018}

@article{Tonkonog:2018aa,
	author = {Dmitry Tonkonog},
	eprint = {1801.06921},
	journal = {arXiv:1801.06921},
	month = {01},
	title = {String topology with gravitational descendants, and periods of {L}andau-{G}inzburg potentials},
	url = {https://arxiv.org/pdf/1801.06921.pdf},
	year = {2018}}

@article{Cieliebak2018,
	author = {Kai Cieliebak and Klaus Mohnke},
	day = {01},
	doi = {10.1007/s00222-017-0767-8},
	issn = {1432-1297},
	journal = {Inventiones mathematicae},
	month = {Apr},
	number = {1},
	pages = {213--295},
	title = {Punctured holomorphic curves and {L}agrangian embeddings},
	url = {https://link.springer.com/content/pdf/10.1007/s00222-017-0767-8.pdf},
	volume = {212},
	year = {2018}}
\Addresses}
\end{document}